\documentclass[aos,nameyear,dvips]{arximspdf}
\usepackage{graphics}
\doi{10.1214/10-AOS792}
\volume{38}
\issue{5}
\pubyear{2010}
\firstpage{2587}
\lastpage{2619}

\makeatletter

\renewcommand{\parallel}{\,\Vert\,}

\newcommand{\dd}{{d}}

\newtheorem{lemma}{Lemma}[section]
\newtheorem{theorem}[lemma]{Theorem}

\newproclaim{example}{Example}

\makeatother

\begin{document}
\begin{frontmatter}

\title{Bayes and empirical-Bayes multiplicity adjustment in the
variable-selection problem}
\runtitle{Bayes and empirical-Bayes multiplicity adjustment}

\begin{aug}
\author[A]{\fnms{James G.} \snm{Scott}\corref{}\thanksref{t1}\ead[label=e1]{james.scott@mccombs.utexas.edu}} and
\author[B]{\fnms{James O.} \snm{Berger}\ead[label=e2]{berger@stat.duke.edu}}
\runauthor{J. G. Scott and J. O. Berger}
\affiliation{University of Texas at Austin and Duke University}
\address[A]{Department of Statistics\\
University of Texas at Austin\\
1 University Station, B6500 \\
Austin, Texas 78712\\
USA \\
\printead{e1}}
\address[B]{Department of Statistics\\
Duke University\\
Box 90251 \\
Durham, North Carolina 27708\\
USA \\
\printead{e2}}
\end{aug}

\thankstext{t1}{Supported in part by the U.S. National Science
Foundation under a Graduate Research Fellowship and
Grants AST-0507481 and DMS-01-03265.}

\received{\smonth{6} \syear{2008}}
\revised{\smonth{11} \syear{2009}}

%
\begin{abstract}
This paper studies the multiplicity-correction effect of standard
Bayesian variable-selection priors in linear regression. Our first
goal is to clarify when, and how, multiplicity correction happens
automatically in Bayesian analysis, and to distinguish this correction
from the Bayesian Ockham's-razor effect. Our second goal is to
contrast empirical-Bayes and fully Bayesian approaches to variable
selection through examples, theoretical results and simulations.
Considerable differences between the two approaches are found. In
particular, we prove a theorem that characterizes a surprising
aymptotic discrepancy between fully Bayes and empirical Bayes. This
discrepancy arises from a different source than the failure to account
for hyperparameter uncertainty in the empirical-Bayes estimate.
Indeed, even at the extreme, when the empirical-Bayes estimate
converges asymptotically to the true variable-inclusion probability,
the potential for a serious difference remains.
\end{abstract}

%

\end{frontmatter}

\section{Introduction}
\label{intro}

This paper addresses concerns about multiplicity in the traditional
variable-selection problem for linear models. We focus on Bayesian and
empirical-Bayesian approaches to the problem. These methods both have
the attractive feature that they can, if set up correctly, account for
multiplicity automatically, without the need for ad-hoc penalties.

Given the huge number of possible predictors in many of today's
scientific problems, these concerns about multiplicity are becoming
ever more relevant. They are especially critical when researchers have
little reason to suspect one model over another, and simply want the
data to flag interesting covariates from a large pool. In such cases,
variable selection is treated less as a formal inferential framework
and more as an exploratory tool used to generate insights about
complex, high-dimensional systems. Still, the results of such studies
are often used to buttress scientific conclusions or guide policy
decisions---conclusions or decisions that may be quite wrong if the
implicit multiple-testing problem is ignored.

Our first objective is to clarify how multiplicity correction enters
Bayesian variable selection: by allowing the choice of prior model
probabilities to depend upon the data in an appropriate way. Some
useful references on this idea include \citet{wallerduncan1969},
\citet{mengdempster1987}, \citet{berry1988},
\citet{westfalljohnsonutts1997}, \citet{berryhochberg1999} and
\citet{scottberger06}. We also clarify the difference between
multiplicity correction and the Bayesian Ockham's-razor effect [see
\citet{jefferysberger92}], which induces a very different type of
penalty on model complexity. This discussion will highlight the fact
that not all Bayesian analyses automatically adjust for multiplicity.

Our second objective is to describe and investigate a peculiar
discrepancy between fully Bayes and empirical-Bayes variable selection.
This discrepancy
seems to arise from a different source than the failure to account for
uncertainty in the empirical-Bayes estimate---the usual issue in such
problems. Indeed, even
when the empirical-Bayes estimate converges asymptotically to the true
hyperparameter value, the potential for a serious difference remains.

The existence of such a discrepancy between fully Bayesian answers and
empirical-Bayes answers---especially one that persists even in the
limit---is of immediate interest to Bayesians, who often use empirical
Bayes as a computational simplification. But the discrepancy is also of
interest to non-Bayesians for at least two reasons.

First, frequentist complete-class theorems suggest that if an
empirical-Bayes analysis does not approximate some fully Bayesian
analysis, then it may be suboptimal and needs alternative
justification. Such justifications can be found for a variety of
situations in \citet{georgefoster2000},
\citet{efrontibshirani2001},
\citet{johnstonesilverman2004},
\citet{bogdanghoshzak2008}, \citet{cuigeorge2006},
\citet{bogdanghosh2008b} and \citet{bogdanetal2008}.

Second, theoretical and numerical investigations of the discrepancy
revealed some unsettling properties of the standard empirical-Bayes
analysis in variable selection. Of most concern is that empirical Bayes
has the potential to collapse to a degenerate solution, resulting in an
inappropriate statement of certainty in the selected regression model.
As a simple example, suppose the usual variable-selection prior is
used, where each variable is presumed to be in the model independently
with an unknown common probability $p$. A common empirical-Bayes method
is to estimate $p$ by marginal maximum likelihood (or Type-II maximum
likelihood, as it is commonly called; see Section \ref{EB-approach}).
This estimated $\hat{p}$ is then used to determine the posterior
probabilities of models. This procedure will be shown to have the
startlingly inappropriate property of assigning final probability 1 to
either the full model or the intercept-only (null) model whenever the
full (or null) model has the largest marginal likelihood, even if this
marginal likelihood is only slightly larger than that of the next-best model.

This is certainly not the first situation in which the Type-II MLE
approach to empirical Bayes has been shown to have problems. But the
unusual character of the problem in variable selection seems not to
have been recognized.

In bringing this issue to light, our goal is not to criticize
empirical-Bayes analysis per se. Indeed, this paper will
highlight many virtues of the empirical-Bayes approach to variable selection,
especially compared to the nonadaptive model prior probabilities that
are often used for variable selection.
Our primary goal is comparative, rather than evaluative, in nature. In
particular, we wish to explore the implications of the above
discrepancy for Bayesians, who are likely to view empirical Bayes as an
approximation to full Bayes analysis, and who wish to understand when
the approximation is a good one. We recognize that others have
alternative goals for empirical Bayes, and that these goals do not
involve approximating full Bayes analysis. Also, there are non-Bayesian
alternatives to marginal maximum likelihood in estimating $p$, as shown
in some of the above papers. The results in this paper suggest that
such alternatives be seriously considered by those wishing to adopt the
empirical-Bayes approach, especially in potentially degenerate situations.

Section \ref{preliminaries} introduces notation. Section \ref
{multiplicity} gives a brief historical and methodological overview of
multiplicity correction for Bayesian variable selection, and focuses on
the issue of clarifying the source and nature of the correction.
Sections~\ref{ebinfogap} and \ref{KLforlinearmodels2} introduce a
theoretical framework for characterizing the differences between fully
Bayesian and empirical-Bayes analyses, and gives several examples and
theoretical results concerning the differences. Section \ref
{numerical-studies} presents numerical results indicating the practical
nature of the differences, through a simulation experiment and a
practical example. Section \ref{summary} gives further discussion of
the results.

\section{Preliminaries}
\label{preliminaries}

\subsection{Notation}

Consider the usual problem of variable selection in linear regression.
Given a vector $\mathbf{Y}$ of $n$ responses and an $n \times m$ design
matrix $\mathbf{X}$, the goal is to select $k$ predictors out of $m$
possible ones for fitting a model of the form
%
%
\begin{equation}
\label{regressionmodel}
Y_i = \alpha+ X_{i j_1} \beta_{j_1} + \cdots+ X_{i j_k} \beta_{j_k} +
\varepsilon_i
\end{equation}
for some $\{j_1, \ldots, j_k \} \subset\{1, \ldots, m\}$, where
$\varepsilon_i \stackrel{\mathrm{i.i.d.}}{\sim} \mathrm{N}(0, \phi^{-1})$
for an unknown
variance $\phi^{-1}$.

All models are assumed to include an intercept term $\alpha$. Let $M_0$
denote the null model with only this intercept term, and let $M_F$
denote the full model with all covariates under consideration. The full
model thus has parameter vector $\bolds{\theta}' = (\alpha,
\bolds{\beta}')$, $\bolds{\beta}' = (\beta_1, \ldots, \beta
_m)'$. Submodels $M_{\bolds{\gamma}}$ are indexed by a binary
vector $\bolds{\gamma}$ of length $m$ indicating a set of
$k_{\bolds{\gamma}} \leq m$ nonzero regression coefficients
$\bolds{\beta}_{\bolds{\gamma}}$:
\[
\gamma_i =
\cases{
0, &\quad if $\beta_i = 0$, \cr
1, &\quad if $\beta_i \neq0$.}
\]

It is most convenient to represent model uncertainty as uncertainty in
$\bolds{\gamma}$, a random variable that takes values in the
discrete space $\{0,1\}^m$, which has $2^m$ members. Inference relies
upon the prior probability of each model, $p(M_{\bolds{\gamma}})$,
along with the marginal likelihood of the data under each model:
%
%
\begin{equation}
\label{marglikemodel}
f(\mathbf{Y} \mid M_{\bolds{\gamma}}) = \int f(\mathbf{Y} \mid
\bolds{\theta}_{\bolds{\gamma}}, \phi)
\pi(\bolds{\theta}_{\bolds{\gamma}},\phi) \,\dd\bolds
{\theta}_{\bolds{\gamma}} \,\dd\phi,
\end{equation}
where $\pi(\bolds{\theta}_{\bolds{\gamma}},\phi)$ is the
prior for model-specific parameters. These together define, up to a
constant, the posterior probability of a model:
%
%
\begin{equation}
\label{postprobmod}
p(M_{\bolds{\gamma}} \mid\mathbf{Y}) \propto p(M_{\bolds
{\gamma}}) f(\mathbf{Y} \mid M_{\bolds{\gamma}}) .
\end{equation}

Let $\mathbf{X}_{\bolds{\gamma}}$ denote the columns of the full
design matrix $\mathbf{X}$ given by the nonzero elements of
$\bolds{\gamma}$, and let $\mathbf{X}^*_{\bolds{\gamma}}$
denote the concatenation $(\mathbf{1}\hspace*{6pt}
\mathbf{X}_{\bolds{\gamma}})$, where $\mathbf{1}$ is a column of
ones corresponding to the intercept $\alpha$. For simplicity, we will
assume that all covariates have been centered so that $\mathbf{1}$ and
$\mathbf{X}_{\bolds{\gamma}}$ are orthogonal. We will also assume
that the common choice $\pi(\alpha)=1$ is made for the parameter
$\alpha$ in each model [see \citet{bergervarsh1998} for a justification
of this choice of prior].

Often all models will have small posterior probability, in which case
more useful summaries of the posterior distribution are quantities such
as the posterior inclusion probabilities of the individual variables:
%
%
\begin{equation}
\label{margincprob}
p_i = \operatorname{Pr}(\gamma_i \neq0 \mid\mathbf{Y}) = \sum
_{\bolds
{\gamma}} 1_{\gamma_i =1} \cdot p(M_{\bolds{\gamma}} \mid\mathbf
{Y}) .
\end{equation}
These quantities also define the median-probability model, which is the
model that includes those covariates having posterior inclusion
probability at least $1/2$. Under many circumstances, this model has
greater predictive power than the most probable model
[\citet{barbieriberger2004}].

\subsection{Priors for model-specific parameters}

There is an extensive body of literature confronting the difficulties
of Bayesian model choice in the face of weak prior information. These
difficulties arise due to the obvious dependence of the marginal
likelihoods in (\ref{marglikemodel}) upon the choice of priors for
model-specific parameters. In general, one cannot use improper priors
on these parameters, since this leaves the resulting Bayes factors
defined only up to an arbitrary multiplicative constant.

This paper chiefly uses null-based $g$-priors [\citet{zellner86}] for
computing the marginal likelihoods in (\ref{marglikemodel}); explicit
expressions can be found in the \hyperref[parameterpriors]{Appendix}. See
\citet{liangpaulo07} for a recent discussion of $g$-priors, and
mixtures thereof, for variable selection.

\section{Approaches to multiple testing}
\label{multiplicity}

\subsection{Bayes factors, Ockham's razor and multiplicity}

In both Bayes and empirical-Bayes variable selection, the marginal
likelihood contains a built-in penalty for model complexity that is
often called the Bayesian ``Ockham's-razor effect''
[\citet{jefferysberger92}]. This penalty arises in integrating the
likelihood across a higher-dimensional parameter space under the more
complex model, resulting in a more diffuse predictive distribution for
the data.

While this is a penalty against more complex models, it is not a
multiple-testing penalty. Observe that the Bayes factor between two
fixed models will not change as more possible variables are thrown into
the mix, and hence will not exert control over the number of false
positives as $m$ grows large.

Instead, multiplicity must be handled through the choice of prior
probabilities of models. The earliest recognition of this idea seems to
be that of Jeffreys in 1939, who gave a variety of suggestions for
apportioning probability across different kinds of model spaces [see
Sections 1.6, 5.0 and 6.0 of \citet{jeffreys1961}, a later
edition]. Jeffreys paid close attention to multiplicity adjustment,
which he called ``correcting for selection.'' In scenarios involving
an infinite sequence of nested models, for example, he recommended
using model probabilities that formed a convergent geometric series, so
that the prior odds ratio for each pair of neighboring models (i.e.,
those differing by a single parameter) was a fixed constant. Another
suggestion, appropriate for more general contexts, was to give all
models of size $k$ a single lump of probability to be apportioned
equally among models of that size. Below, in fact, the fully Bayesian
solution to multiplicity correction will be shown to have exactly this
flavor.

It is interesting that, in the variable-selection problem, assigning
all models equal prior probability (which is equivalent to assigning
each variable prior probability of $1/2$ of being in the model) provides
no multiplicity control. This is most obvious in the orthogonal
situation, which can be viewed as $m$ independent tests of $H_i\dvtx
\beta_i=0$. If each of these tests has prior probability of $1/2$, there
will be no multiplicity control as $m$ grows. Indeed, note that this
``pseudo-objective'' prior reflects an a priori expected model
size of $m/2$ with a standard deviation of $\sqrt{m}/2$, meaning that
the prior for the fraction of included covariates becomes very tight
around $1/2$ as $m$ grows. See \citet{bogdanetal2008} for
extensive discussion of this issue.

\subsection{Variable-selection priors and empirical Bayes}
\label{EB-approach}

The standard modern practice in Bayesian variable-selection problems is
to treat variable inclusions as exchangeable Bernoulli trials with
common success probability $p$, which implies that the prior
probability of a model is given by
%
%
\begin{equation}
\label{conditionalmodelprior}
p(M_{\bolds{\gamma}} \mid p) = p^{k_{\bolds{\gamma}}}
(1-p)^{m - k_{\bolds{\gamma}}}
\end{equation}
with $k_{\bolds{\gamma}}$ representing the number of included
variables in the model.

We saw above that selecting $p=1/2$ does not provide multiplicity
correction. Treating $p$ as an unknown parameter to be estimated from
the data will, however, yield an automatic multiple-testing penalty.
The intuition is that, as $m$ grows with the true $k$ remaining fixed,
the posterior distribution of $p$ will concentrate near $0$, so that
the situation is the same as if one had started with a very low prior
probability that a variable should be in the model
[\citet{scottberger06}]. Note that one could adjust for multiplicity
subjectively, by specifying $p$ to reflect subjective belief in the
proportion of variables that should be included. No fixed choice of $p$
that is independent of $m$, however, can adjust for multiplicity.

The empirical-Bayes approach to variable selection was popularized by
\citet{georgefoster2000}, and is a common strategy for treating
the prior inclusion probability $p$ in (\ref{conditionalmodelprior}) in
a data-dependent way. The most common approach is to estimate the prior
inclusion probability by maximum likelihood, maximizing the marginal
likelihood of $p$ summed over model space (often called Type-II maximum
likelihood):
%
%
\begin{equation}
\label{EBmaximization}
\hat{p} = \mathop{\arg\max}_{p \in[0,1]} \sum_{\bolds{\gamma}} p(
M_{\bolds{\gamma}} \mid p) \cdot f(\mathbf{Y} \mid M_{\bolds
{\gamma}}) .
\end{equation}
One uses this in (\ref{conditionalmodelprior}) to define the
ex-post prior probabilities $p(M_{\bolds{\gamma}} \mid\hat
{p}) =
\hat{p}^{k_{\bolds{\gamma}}} (1 - \hat{p})^{m - k_{\bolds
{\gamma}}}$, resulting in final model posterior probabilities
%
%
\begin{equation}
\label{eb-posteriors}
p(M_{\bolds{\gamma}} \mid\mathbf{Y}) \propto\hat
{p}^{k_{\bolds{\gamma}}} \cdot(1 - \hat{p})^{m - k_{\bolds
{\gamma}}} f(\mathbf{Y} \mid M_{\bolds{\gamma}}) .
\end{equation}
The EB solution $\hat{p}$ can be found either by direct numerical
optimization or by the EM algorithm detailed in
\citet{liangpaulo07}. For an overview of empirical-Bayes
methodology, see \citet{carlinlouis2000}.


It is clear that the empirical-Bayes approach will control for
multiplicity in a straightforward way: if there are only $k$ true
variables and $m$ grows large, then $\hat{p} \rightarrow0$. This will
make it increasingly more difficult for all variables to overcome the
ever-stronger prior bias against their relevance.

\subsection{A fully Bayesian version}

Fully Bayesian variable-selection priors have been discussed by
\citet{leysteel2007}, \citet{cuigeorge2006} and
\citet{carvalhoscott2007}, among others. These priors assume that
$p$ has a Beta distribution, $p \sim\operatorname{Be}(a, b)$, giving
%
%
\begin{equation}
\label{marginalmodelprior}
p(M_{\bolds{\gamma}}) = \int_0^1 p(M_{\bolds{\gamma}} \mid p)
\pi(p) \,\dd p = \frac{\beta(a + k_{\bolds{\gamma}}, b + m -
k_{\bolds{\gamma}})}{\beta(a,b)} ,
\end{equation}
where $\beta(\cdot,\cdot)$ is the beta function. For the default choice
of $a = b = 1$, implying a uniform prior on $p$, this reduces to
%
%
\begin{equation}
\label{marginalmodelpriorunif}
p(M_{\bolds{\gamma}}) = \frac{(k_{\bolds{\gamma}})! (m -
k_{\bolds{\gamma}})!}{(m+1)(m!)}= \frac{1}{m+1}\pmatrix{m \cr
k_{\bolds{\gamma}}}^{-1} .
\end{equation}

%
%
\begin{figure}

\includegraphics{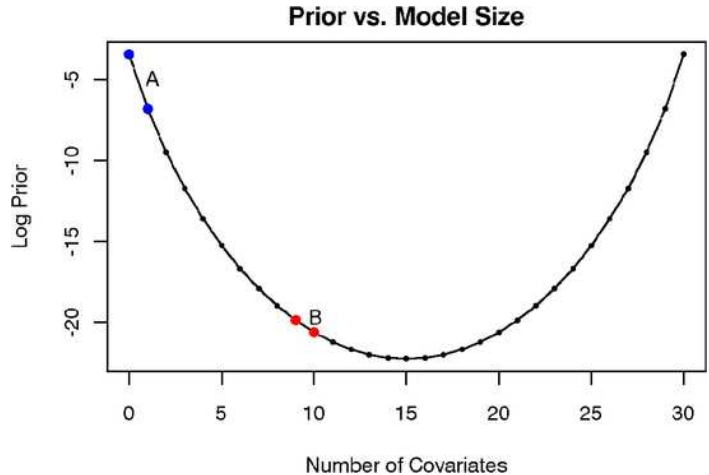}

\caption{Prior model probability versus model size.}\label{Uplot}
\end{figure}

%
%
\begin{figure}[b]

\includegraphics{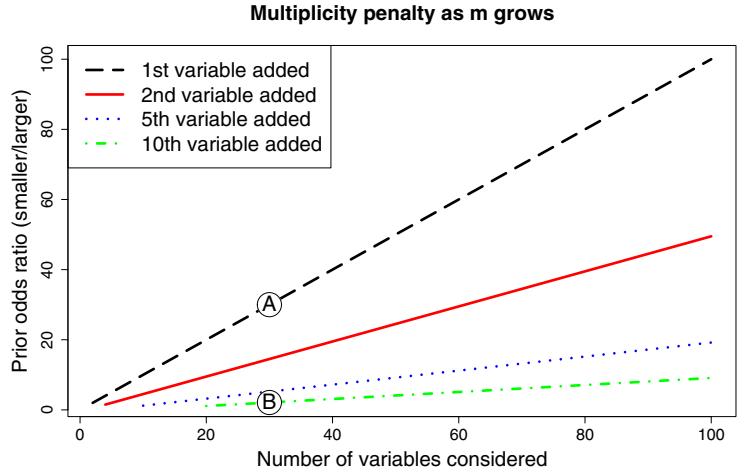}

\caption{Multiplicity penalties as $m$ grows.}\label{penaltyplot}
\end{figure}

We call these expressions deriving from the uniform prior on $p$ the
``fully Bayes'' version of variable selection priors, though of course
many other priors could be used (including those incorporating
subject-area information). Utilizing these prior probabilities in (\ref
{postprobmod}) yields the following posterior probabilities:
%
%
\begin{equation}
\label{postprobmodfull}
p(M_{\bolds{\gamma}} \mid\mathbf{Y}) \propto\frac{1}{m+1}\pmatrix{m
\cr k_{\bolds{\gamma}}}^{-1} f(\mathbf{Y} \mid M_{\bolds
{\gamma}}) .
\end{equation}
This has the air of paradox: in contrast to (\ref{eb-posteriors}),
where the multiplicity adjustment is apparent, here $p$ has been
marginalized away. How can $p$ then be adjusted by the data so as to
induce a multiplicity-correction effect?

Figures \ref{Uplot} and \ref{penaltyplot} hint at the answer, which is
that the multiplicity penalty was always in the prior probabilities in
(\ref{marginalmodelpriorunif}) to begin with; it was just hidden. In
Figure~\ref{Uplot} the prior log-probability is plotted as a function
of model size for a particular value of $m$ (in this case $30$). This
highlights the marginal penalty that one must pay for adding an extra
variable: in moving from the null model to a model with one variable,
the fully Bayesian prior favors the simpler model by a factor of $30$
(label A). This penalty is not uniform: models of size 9, for example,
are favored to those of size 10 by a factor of only $2.1$ (label B).

Figure \ref{penaltyplot} then shows these penalties getting steeper as
one considers more models. Adding the first variable incurs a
$30$-to-$1$ prior-odds penalty if one tests 30 variables (label A as
before), but a $60$-to-$1$ penalty if one tests 60 variables.
Similarly, the 10th-variable marginal penalty is about two-to-one for
30 variables considered (label B), but would be about four-to-one for
60 variables.

We were careful above to distinguish this effect from the
Ockham's-razor penalty coming from the marginal likelihoods. But
marginal likelihoods are clearly relevant. They determine where models
will sit along the curve in Figure \ref{Uplot}, and thus will determine
whether the prior-odds multiplicity penalty for adding another variable
to a good model will be more like $2$, more like $30$ or something else
entirely. Indeed, note that, if only large models have significant
marginal likelihoods, then the ``multiplicity penalty'' will now become a
``multiplicity advantage,'' as one is on the increasing part of the curve
in Figure \ref{Uplot}. (This is also consistent with the
empirical-Bayes answer: if $\hat p > 0.5$, then the analysis will
increase the chance of variables entering the model.)

Interestingly, the uniform prior on $p$ also gives every variable a
marginal prior inclusion probability of $1/2$; these marginal
probabilities are the same as those induced by the ``pseudo-objective''
choice of $p = 1/2$. Yet because probability is apportioned among
models in a very different way, profoundly different behaviors emerge.

%
%
\begin{table}
\caption{Posterior inclusion probabilities
($\times$100) for the 10 real variables in the simulated data, along
with the number of false positives (posterior inclusion probability
greater than $1/2$) from the ``pure noise'' columns in the design matrix.
Marginal likelihoods were calculated (under Zellner--Siow priors) by
enumerating the model space in the $m=11$ and $m=20$ cases, and by 5
million iterations of the feature-inclusion stochastic-search algorithm
[Berger and Molina (\protect\citeyear{bergermolina2005}),
Scott and Carvalho (\protect\citeyear{scottcarvalho2007b})] in the $m=50$ and $m=100$
cases}\label{postincfakedata}
\begin{tabular*}{\tablewidth}{@{\extracolsep{\fill}}lcccccccccccccccc@{}}
\hline
& \multicolumn{16}{c@{}}{\textbf{Method and number of noise variables}}
\\[-4pt]
& \multicolumn{16}{c@{}}{\hrulefill} \\
&\multicolumn{4}{c}{\textbf{Uncorrected}} & \multicolumn{4}{c}{\textbf{Fully Bayes}}
& \multicolumn{4}{c}{\textbf{Oracle Bayes}} & \multicolumn{4}{c@{}}{\textbf{Empirical
Bayes}} \\[-4pt]
&\multicolumn{4}{c}{\hrulefill} & \multicolumn{4}{c}{\hrulefill}
& \multicolumn{4}{c}{\hrulefill} & \multicolumn{4}{c@{}}{\hrulefill}\\
\multicolumn{1}{@{}l}{\textbf{Signal}} & \textbf{1} & \textbf{10}
& \textbf{40} & \textbf{90} & \textbf{1} & \textbf{10} & \textbf{40}
& \textbf{90} & \textbf{1} & \textbf{10} & \textbf{40} & \textbf{90}
& \textbf{1} & \textbf{10} & \textbf{40} & \textbf{90}\\
\hline
$-1.08$ & 99 & 99 & 99 & 99 & 99 & 99 & 99 & 99 & 99 & 99 & 99 & 99 &
99 & 99 & 99 & 99 \\
$-0.84$ & 99 & 99 & 99 & 99 & 99 & 99 & 99 & 98 & 99 & 99 & 99 & 99 &
99 & 99 & 99 & 99 \\
$-0.74$ & 99 & 99 & 99 & 99 & 99 & 99 & 99 & 99 & 99 & 99 & 99 & 99 &
99 & 99 & 99 & 99\\
$+0.63$ & 99 & 99 & 99 & 99 & 99 & 99 & 92 & 73 & 99 & 99 & 97 & 87 &
99 & 99 & 93 & 80\\
$-0.51$ & 97 & 97 & 99 & 99 & 91 & 94 & 71 & 34 & 99 & 97 & 85 & 52 &
93 & 95 & 74 & 44\\
$+0.41$ & 92 & 91 & 99 & 99 & 96 & 86 & 56 & 22 & 99 & 91 & 72 & 35 &
97 & 88 & 60 & 25 \\
$+0.35$ & 77 & 77 & 99 & 99 & 89 & 68 & 30 & 05 & 97 & 77 & 45 & 11 &
91 & 72 & 35 & 07 \\
$-0.30$ & 29 & 28 & 28 & 12 & 55 & 24 & 04 & 00 & 79 & 28 & 06 & 01 &
64 & 25 & 04 & 01 \\
$+0.18$ & 26 & 28 & 24 & 27 & 51 & 25 & 03 & 01 & 79 & 28 & 04 & 01 &
62 & 24 & 04 & 01 \\
$+0.07$ & 21 & 24 & 05 & 01 & 45 & 21 & 03 & 01 & 70 & 24 & 05 & 01 &
56 & 22 & 03 & 01 \\
[3pt]
FPs & 0 & 2 & 5 & 10 & 0 & 1 & 0 & 0 & 0 & 2 & 1 & 0 & 0 & 1 & 1 & 0 \\
\hline
\end{tabular*}
\end{table}

For example, Table \ref{postincfakedata} compares these two regimes on
a simulated data set for which the true value of $k$ was fixed at $10$.
The goal of the study is, in essence, to understand how posterior
probabilities adapt to situations of increasingly egregious ``data
dredging,'' where a set of true covariates is tested in the presence of
an ever-larger group of spurious covariates. We used a simulated
$m=100$ design matrix of $\mathrm{N}(0,1)$ covariates and $10$ regression
coefficients that differed from zero, along with $90$ coefficients that
were identically zero. The table summarizes the posterior inclusion
probabilities of the $10$ real variables as we test them along with an
increasing number of noise variables (first 1, then 10, 40 and 90). It
also indicates how many false positives (defined as having posterior
inclusion probability $\geq0.5$) are found among the noise variables.
Here, ``uncorrected'' refers to giving all models equal prior
probability by setting $p = 1/2$. ``Oracle Bayes'' is the result from
choosing $p$ to reflect the known fraction of nonzero covariates.

The following points can be observed:
\begin{itemize}
\item The fully Bayes and empirical Bayes procedures both exhibit clear
multiplicity adjustment: as the number of noise variables increases,
the posterior inclusion probabilities of variables decrease. The
uncorrected Bayesian analysis shows no such adjustment and can, rather
bizarrely, sometimes have the posterior inclusion probabilities
increase as noise variables are added.
\item On the simulated data, proper multiplicity adjustment yields
reasonably strong control over false positives, in the sense that the
number of false positives appears bounded (and small) as $m$ increases.
In contrast, the number of false positives appears to be increasing
linearly for the uncorrected Bayesian analysis, as would be expected.
\item The full Bayes, empirical Bayes and oracle Bayes answers are all
qualitatively (though not quantitatively) similar; indeed, if one
adopted the (median probability model) prescription of selecting those
variables with posterior inclusion probability greater than $1/2$, they
would both always select the same variables, except in two instances.
\end{itemize}


%
%
\begin{table}
\tablewidth=276pt
\caption{Posterior inclusion probabilities for
the important main effects, quadratic effects and cross-product effects
for ozone-concentration data under $g$-priors. Key: $p=1/2$ implies
that all models have equal prior probability; FB is fully Bayes; EB is
empirical Bayes}\label{ozoneinclusion}
\begin{tabular*}{\tablewidth}{@{\extracolsep{\fill}}lccc@{}}
\hline
& {$\bolds{p=1/2}$} & \textbf{FB} & \textbf{EB} \\
\hline
x1 & 0.83 & 0.42 & 0.54\\
x2 & 0.13 & 0.03 & 0.05\\
x3 & 0.09 & 0.02 & 0.03\\
x4 & 0.94 & 0.73 & 0.84\\
x5 & 0.33 & 0.06 & 0.10\\
x6 & 0.38 & 0.07 & 0.10\\
x7 & 0.34 & 0.36 & 0.29\\
x8 & 0.78 & 0.74 & 0.77\\
x9 & 0.20 & 0.03 & 0.05\\
x10 & 0.96 & 0.96 & 0.97\\
x1--x1 & 1.00 & 0.97 & 0.99\\
x9--x9 & 0.95 & 0.82 & 0.91\\
x1--x2 & 0.48 & 0.16 & 0.24\\
x4--x7 & 0.33 & 0.10 & 0.15\\
x6--x8 & 0.43 & 0.25 & 0.34\\
x7--x8 & 0.31 & 0.13 & 0.18\\
x7--x10 & 0.71 & 0.86 & 0.85\\
\hline
\end{tabular*}
\end{table}

The differences between corrected and uncorrected analyses are quite
stark, and calls into question the use of nonadaptive priors in
situations with large numbers of potentially spurious covariates. For
example, Table \ref{ozoneinclusion} shows the posterior inclusion
probabilities for a model of ozone concentration levels outside Los
Angeles that includes 10 atmospheric variables along with all squared
terms and second-order interactions ($m = 65$). Probabilities are given
for uncorrected ($p=1/2$), empirical Bayes and fully Bayesian analyses.
All variables appear uniformly less impressive when adjusted for multiplicity.

Other examples of such multiplicity correction put into practice can be
found throughout the literature. For nonparametric problems, see
\citet{gopalanberry1998}; for gene-expression studies, see
\citet{domuller2005}; for econometrics, see
\citet{leysteel2007}; for Gaussian graphical models, see
\citet{carvalhoscott2007}; and for time-series data, see
\citet{scott2008}.

\section{Theoretical comparison of Bayes and empirical Bayes}
\label{ebinfogap}

\subsection{Motivation}

The previous section showed some examples where fully Bayes and
empirical-Bayes methods gave qualitatively similar results. While this
rough correspondence between the two approaches does seem to hold in a
wide variety of applied problems, we now turn attention to the question
of when, and how, it fails.

We begin with a surprising lemma that indicates the need for caution
with empirical-Bayes methods in variable selection. The lemma refers to
the variable-selection problem, with the prior variable inclusion
probability $p$ being estimated by marginal (or Type-II) maximum
likelihood in the empirical-Bayes approach.
\begin{lemma} \label{EBextremelemma}
In the variable-selection problem, if $M_0$ has the (strictly) largest
marginal likelihood, then the Type-\textup{II} MLE estimate of $p$ is $\hat{p} =
0$. Similarly, if $M_F$ has the (strictly) largest marginal likelihood,
then $\hat{p} = 1$.
\end{lemma}
\begin{pf}
Since $p(M_{\bolds{\gamma}})$ sums to $1$ over $\bolds{\gamma
}$, the marginal likelihood of $p$ satisfies
%
%
\begin{equation}
\label{maximizing}
f(\mathbf{Y}) = \sum_{\Gamma} f(\mathbf{Y} \mid M_{\bolds{\gamma
}}) p(M_{\bolds{\gamma}})
\leq
\max_{\bolds{\gamma} \in\Gamma} f(\mathbf{Y} \mid M_{\bolds
{\gamma}}) .
\end{equation}
Furthermore, the inequality is strict under the conditions of the lemma
(because the designated marginals are strictly largest), unless the
prior assigns $p(M_{\bolds{\gamma}}) = 1$ to the maximizing
marginal likelihood. The only way that $p(M_{\bolds{\gamma}}) =
p^{k_{\bolds{\gamma}}} \cdot(1-p)^{m - k_{\bolds{\gamma}}}$
can equal 1 is for $p$ to be 0 or 1 and for the model to be $M_0$ or
$M_F$, respectively. At these values of $p$, equality is indeed
achieved in (\ref{maximizing}) under the stated conditions, and the
results follow.
\end{pf}

As a consequence, the empirical-Bayes approach here would assign final
probability 1 to $M_0$ whenever it has the largest marginal likelihood,
and final probability~1 to $M_F$ whenever it has the largest marginal
likelihood. These are clearly very unsatisfactory answers.

The above lemma does highlight a specific, undesirable property of the
empirical-Bayes approach to variable selection---one whose practical
significance we investigate by simulation in Section \ref
{numerical-studies}. For the most part, however, the rest of our
results are of a fundamentally different character. We will not be
evaluating either the fully Bayes or the empirical-Bayes approach
according to an objective yardstick, such as how well each one does at
recovering true relationships or suppressing false ones. Instead, we
focus on comparing the two approaches to each other in a more formal
way. As mentioned above, our fundamental goal is to understand when,
and how, empirical Bayes corresponds asymptotically to full Bayes
analysis. Such a comparison is certainly of interest, both to Bayesians
who might consider empirical Bayes as a computational approximation,
and to frequentists for the reasons mentioned in the \hyperref[intro]{Introduction}.

To explore the difference between these two approaches, it is useful to
abstract the problem somewhat and suppose simply that the data $\mathbf
{Y}$ have sampling density $f(\mathbf{Y} \mid\bolds{\theta})$,
and let $\bolds{\theta} \in\Theta$ have prior density $\pi
(\bolds{\theta} \mid\bolds{\lambda})$ for some unknown
hyperparameter $\bolds{\lambda} \in\Lambda$. Empirical-Bayes
methodology typically proceeds by estimating $\bolds{\lambda}$
from the data using a consistent estimator. [The Type-II MLE approach
would estimate $\lambda$ by the maximizer of the marginal likelihood
$m(\mathbf{Y} \mid\bolds{\lambda}) = \int_{\Lambda} f(\mathbf{Y}
\mid\bolds{\theta})\pi(\bolds{\theta} \mid\bolds
{\lambda}) \,d \bolds{\theta}$, and this will typically be
consistent in empirical-Bayes settings.] It is then argued that (at
least asymptotically) the Bayesian analysis with $\hat
{\bolds{\lambda}}$ will be equivalent to the Bayesian analysis\vspace*{1pt} if one knew
$\bolds{\lambda}$. (This claim is most interesting when the prior
for $\hat{\bolds{\lambda}}$ is unknown; if it is known, then there
are also strong frequentist reasons to use this prior in lieu of
empirical Bayes.)

To contrast this with a full Bayesian analysis, suppose we have a prior
density $\pi(\bolds{\lambda})$ for $\bolds{\lambda}$ and a
target function $\psi(\bolds{\theta}, \mathbf{Y} \mid\bolds
{\lambda})$. For instance, $\psi$ could be the posterior mean of
$\bolds{\theta}$ given $\bolds{\lambda}$ and $\mathbf{Y}$, or
it could be the conditional posterior distribution of $\bolds
{\theta}$ given $\bolds{\lambda}$ and $\mathbf{Y}$. The
empirical-Bayesian claim, in this context, would be that
%
%
\begin{equation}
\label{EBpost}
\int_{\Lambda} \psi(\bolds{\theta}, \mathbf{Y} \mid\bolds
{\lambda}) \pi(\bolds{\lambda} \mid\mathbf{Y}) \,\dd\bolds
{\lambda}
\approx
\psi(\bolds{\theta}, \mathbf{Y} \mid\hat{\bolds{\lambda}}) ,
\end{equation}
that is, that the full Bayesian answer on the left can be well
approximated by the empirical-Bayes answer on the right.
The justification for (\ref{EBpost}) would be based on the fact that,
typically, $\pi(\bolds{\lambda} \mid\mathbf{Y})$ will collapse to
a point mass near the true $\bolds{\lambda}$ as the sample size
increases, so that (\ref{EBpost}) will hold for appropriately smooth
functions $\psi(\bolds{\theta}, \mathbf{Y} \mid\bolds{\lambda
})$ when the sample size is large.

There are typically better approximations to the left-hand side of
(\ref
{EBpost}), such as the Laplace approximation. These, however, are
focused on reproducing the full-Bayes analysis through an analytic
approximation, and are not ``empirical-Bayes'' per se. Likewise,
higher-order empirical-Bayes analysis will likely yield better results here,
but the issue is in realizing when one needs to resort to such
higher-order analysis in the first place, and in understanding why this
is so for problems such as variable selection.

That (\ref{EBpost}) could fail for nonsmooth $\psi(\bolds{\theta},
\mathbf{Y} \mid\bolds{\lambda})$ is no surprise. But what may
come as a surprise is that this failure can also occur for very common
functions.

Most notably, it fails for the conditional posterior density itself.
Indeed, in choosing $\psi(\bolds{\theta}, \mathbf{Y} \mid
\bolds{\lambda}) = \pi(\bolds{\theta} \mid\bolds{\lambda
}, \mathbf{Y})$, the left-hand side of (\ref{EBpost}) is just the
posterior density of $\bolds{\theta}$ given $\mathbf{Y}$, which
(by definition) can be written as
%
%
\begin{equation}
\label{Fullpost}
\pi_F(\bolds{\theta} \mid\mathbf{Y}) \propto f(\mathbf{Y} \mid
\bolds{\theta}) \int_{\Lambda} \pi(\bolds{\theta} \mid
\bolds{\lambda}) \pi(\bolds{\lambda}) \,\dd\bolds
{\lambda} .
\end{equation}
On the other hand, for this choice of $\psi$, (\ref{EBpost}) becomes
%
%
\begin{equation}
\label{EBpost2}
\pi_E(\bolds{\theta} \mid\mathbf{Y}) \approx\pi(\bolds
{\theta} \mid\mathbf{Y}, \hat{\bolds{\lambda}}) \propto f(\mathbf
{Y} \mid\bolds{\theta}) \cdot\pi(\bolds{\theta} \mid\hat
{\bolds{\lambda}}) ,
\end{equation}
and the two expressions on the right-hand sides of (\ref{Fullpost}) and
(\ref{EBpost2}) can be very different. [This difference may not matter,
of course; for instance, if $f(\mathbf{Y} \mid\bolds{\theta})$ is
extremely concentrated as a likelihood, the prior
being used may not matter.]

As an indication as to what goes wrong in (\ref{EBpost}) for this
choice of $\psi$, note that
%
%
\begin{eqnarray}
\pi_F(\bolds{\theta} \mid\mathbf{Y})
&=& \int_{\Lambda} \pi(\bolds{\theta} \mid\bolds{\lambda} ,
\mathbf{Y}) \cdot\pi(\bolds{\lambda} \mid\mathbf{Y})
\,\dd\bolds{\lambda} \nonumber\\
&=& \int_{\Lambda} \frac{\pi(\bolds{\theta}, \bolds{\lambda}
\mid\mathbf{Y})}{\pi(\bolds{\lambda} \mid\mathbf{Y})} \cdot
\pi(\bolds{\lambda} \mid\mathbf{Y}) \,\dd\bolds{\lambda} \\
&=& \int_{\Lambda} \frac{f(\mathbf{Y} \mid\bolds{\theta}) \pi
(\bolds{\theta} \mid\bolds{\lambda}) \pi(\bolds{\lambda
})}{f(\mathbf{Y})\pi(\bolds{\lambda} \mid\mathbf{Y})} \cdot
\pi(\bolds{\lambda} \mid\mathbf{Y}) \,\dd\bolds{\lambda} .
\end{eqnarray}
Of course, these elementary calculations simply lead to (\ref
{Fullpost}) after further algebra. But they illuminate the fact that,
while $\pi(\bolds{\lambda} \mid\mathbf{Y})$ may indeed be
collapsing to a point mass at the true $\bolds{\lambda}$, this
term occurs in both the numerator and the denominator of the integrand
and therefore cancels. The accuracy with which a point mass at $\hat
{\bolds{\lambda}}$ approximates $\pi(\bolds{\lambda} \mid
\mathbf{Y})$ is thus essentially irrelevant from the standpoint of full
Bayes analysis.

\subsection{Comparison using Kullback--Leibler convergence}

Our goal, then, is to understand when (\ref{Fullpost}) and (\ref
{EBpost2}) will yield the same answers, in an asymptotic sense. The
closeness of these two distributions will be measured by
Kullback--Leibler divergence, a standard measure for comparing a pair
of distributions $P$ and $Q$ over parameter space $\Theta$:
%
%
\begin{equation}
\operatorname{KL}(P \parallel Q) = \int_{\Theta} P(\bolds{\theta})
\log
\biggl( \frac{P(\bolds{\theta})}{Q(\bolds{\theta})} \biggr)
\,\dd\bolds{\theta} .
\end{equation}

Kullback--Leibler divergence can be used to formalize the notion of
empirical-Bayes convergence to fully Bayesian analysis as follows:

\subsubsection*{KL empirical-Bayes convergence}
Suppose the data $\mathbf{Y}$ and parameter $\bolds{\theta}$ have
joint distribution $p(\mathbf{Y}, \bolds{\theta} \mid\bolds
{\lambda})$, where $\bolds{\theta} \in\Theta$ is of dimension
$m$, and where $\bolds{\lambda} \in\Lambda$ is of fixed dimension
that does not grow with $m$. Let $\pi_E = \pi(\psi(\bolds{\theta})
\mid\mathbf{Y}, \hat{\bolds{\lambda}})$ be the empirical-Bayes
posterior distribution for some function of the parameter $\psi
(\bolds{\theta})$, and let $\pi_F = \pi(\psi(\bolds{\theta})
\mid\mathbf{Y}) = \int_{\Lambda} \pi(\psi(\bolds{\theta}) \mid
\mathbf{Y}, \bolds{\lambda}) \cdot\pi(\bolds{\lambda}) \,\dd
\bolds{\lambda}$ be the corresponding fully Bayesian posterior
under the prior $\pi(\bolds{\lambda})$. If, for every $\bolds
{\lambda} \in\Lambda$, $\operatorname{KL}(\pi_F \parallel\pi_E)
\to0$ in
probability [expectation] under $p(\mathbf{Y}, \bolds{\theta} \mid
\bolds{\lambda})$ as $m \to\infty$, then $\pi_E$ will be said to
be KL-convergent in probability [expectation] to the fully Bayesian
posterior $\pi_F$.

Note that KL convergence is defined with respect to a particular
function of the parameter, along with a particular prior distribution
on the hyperparameter. The intuition is the following. Suppose that in
trying to estimate a given function $\psi(\bolds{\theta})$, it is
possible to construct a reasonable prior $\pi(\bolds{\lambda})$
such that the KL-convergence criterion is met. Then the empirical Bayes
and full Bayes analysis will disagree for every finite sample size, but
are at least tending toward agreement asymptotically. If, on the other
hand, it is not possible to find a reasonable prior $\pi(\bolds
{\lambda})$ that leads to KL convergence, then estimating $\psi
(\bolds{\theta})$ by empirical Bayes is dubious from the fully
Bayesian perspective. A Bayesian could not replicate such a procedure
even asymptotically, while a frequentist may be concerned by
complete-class theorems. (A ``reasonable'' prior is a necessarily vague
notion, but obviously excludes things such as placing a point mass at
$\hat{\bolds{\lambda}}$.)

Instead of KL divergence, of course, one might instead use another
distance or divergence measure. The squared Hellinger distance is one
such possibility:
\[
\mathrm{H}^2(P \parallel Q) = \frac{1}{2} \int_{\Theta} \bigl( \sqrt
{P(\bolds{\theta})} - \sqrt{Q(\bolds{\theta})} \bigr)^2
\,\dd\bolds{\theta} .
\]
Most of the subsequent results, however, use KL divergence because of
its familiarity and analytical tractability.

\subsection{An orthogonal example}
\label{normalmeanssection}

As a simple illustration of the above ideas, consider the following two
examples of empirical-Bayes analysis. The first example satisfies the
convergence criterion; the second does not. Both examples concern the
same sampling model, in which we observe a series of conditionally
independent random variables $y_i \sim\mathrm{N}(\theta_i, 1)$, and
where we
know that $\theta_i \sim\mathrm{N}(\mu, 1)$. Thus, the hyperparameter
$\lambda
=\mu$ here. Let $\bolds{\theta} = (\theta_1, \ldots,\theta_m)$ and
$\mathbf{y} = (y_1, \ldots, y_m)$.

Alternatively, this can be thought of as an orthogonal regression
problem where both the dimension and number of samples are growing at
the same rate: $\mathbf{y}= X \bolds{\theta}+ \varepsilon$, with
$X$ being the $m
\times m$ identity matrix. This framing makes the connection to
variable selection much more plain.

The natural empirical-Bayes estimate of $\mu$ is the sample mean
$\hat{\mu}_E = \bar{y}$, which is clearly consistent for $\mu$ as $m
\to\infty$ and converges at the usual $1/\sqrt{m}$ rate. A~standard
hyperprior in a fully Bayesian analysis, on the other hand, would be
$\mu\sim\mathrm{N}(0,A)$ for some specified $A$; the objective hyperprior
$\pi(\mu) =1$ is essentially the limit of this as $A \rightarrow
\infty$. Using the expressions given in, for example,
\citet{bergerbook2ed}, the empirical-Bayes and full Bayes
posteriors are
%
%
\begin{eqnarray}\qquad
\label{generalEBpost}
\pi_E(\bolds{\theta} \mid\mathbf{y}, \hat{\mu}_E) &=& \mathrm{N}\bigl(
\tfrac{1}{2} (\mathbf{y} + \bar{y} \mathbf{1} ), \tfrac
{1}{2} \mathbf{I} \bigr), \\
\label{generalFBpost}
\pi_F(\bolds{\theta} \mid\mathbf{y}) &=&
\mathrm{N}
\biggl( \frac{1}{2} (\mathbf{y} + \bar{y} \mathbf{1} ) -
\biggl( \frac{1}{mA + 2} \biggr) \bar{y} \mathbf{1}, \frac{1}{2 } \mathbf
{I} + \frac{A}{2(mA + 2)} (\mathbf{1} \mathbf{1}^t) \biggr) ,
\end{eqnarray}
where $\mathbf{I}$ is the identity matrix and $\mathbf{1}$ is a column
vector of all ones.
\begin{example}\label{Example1} Suppose only the first normal mean,
$\theta
_1$, is of interest, meaning that the target function $\psi(\bolds
{\theta}) = \theta_1$. Then sending $A \to\infty$ yields
%
%
\begin{eqnarray}
\pi_E(\theta_1 \mid\mathbf{y}, \hat{\mu}_E) &=& \mathrm{N}( [y_1
+ \bar
{y}]/2, 1/2 ), \\
\pi_F(\theta_1 \mid\mathbf{y}) &=& \mathrm{N}( [y_1 + \bar{y}]/2,
1/2 +
[2m]^{-1} ) .
\end{eqnarray}
It is easy to check that $\operatorname{KL}(\pi_F \parallel\pi_E)
\to0$ as $m
\rightarrow\infty$. Hence, $\pi_E(\theta_1)$ arises from a
KL-convergent EB procedure under a reasonable prior, since it
corresponds asymptotically to the posterior given by the objective
prior on the hyperparameter~$\mu$.
\end{example}
\begin{example}\label{Example2} Suppose now that $\bolds{\theta}$, the
entire vector of means, is of interest [hence, $\psi(\bolds{\theta
}) = \bolds{\theta}$]. The relevant distributions are then the
full $\pi_E$ and $\pi_F$ given in (\ref{generalEBpost}) and (\ref
{generalFBpost}), with parameters $(\hat{\bolds{\theta}}_E, \Sigma
_E)$ and $(\hat{\bolds{\theta}}_F, \Sigma_F)$, respectively.

A straightforward computation shows that $\operatorname{KL}(\pi_F
\parallel\pi
_E)$ is given by
%
%
\begin{eqnarray} \label{normalmeansKL}
\operatorname{KL} &=& \frac{1}{2} \biggl[ \log\biggl( \frac{\det\Sigma_E}{\det
\Sigma_F} \biggr) + \operatorname{tr} (\Sigma_E^{-1} \Sigma_F) + (\hat
{\bolds{\theta}}_E - \hat{\bolds{\theta}}_F)^t \Sigma
_E^{-1} (\hat{\bolds{\theta}}_E - \hat{\bolds{\theta}}_F) -
m \biggr] \hspace*{-22pt}\nonumber\\[-8pt]\\[-8pt]
&=& \frac{1}{2} \biggl[ - \log\biggl( 1 + \frac{mA}{mA + 2} \biggr) +
\frac{mA}{mA + 2} + 2m \biggl( \frac{1}{mA + 2} \biggr)^2 \bar{y}^2
\biggr].\nonumber
\end{eqnarray}
For any nonzero choice of $A$ and for any finite value of the
hyperparameter $\mu$, it is clear that under $p(\mathbf{y}, \bolds
{\theta} \mid\mu)$ the quantity $[2m / (mA + 2)^2] \cdot\bar{y}^2
\to
0$ in probability as $m \to\infty$. Hence, for any value of $A$
(including $A = \infty$), the KL divergence in (\ref{normalmeansKL})
converges to $(1 - \log2)/2 > 0$.

Of course, this only considers priors of the form $\mu\sim\mathrm{N}(0,A)$,
but the asymptotic normality of the posterior for $\mu$ can be used to
prove the result for essentially any prior that satisfies the usual
regularity conditions, suggesting that there is no reasonable prior for
which $\pi_E(\bolds{\theta})$ is KL-convergent.
\end{example}

The crucial difference here is that, in the second example, the
parameter of interest increases in dimension as information about the
hyperparameter $\mu$ accumulates. This is not the usual situation in
asymptotic analysis. Hence, even as $\hat{\bolds{\theta}}_E$ and
$\hat{\bolds{\theta}}_F$ are getting closer to each other
elementwise, the KL divergence does not shrink to 0 as expected.

Two further comments are in order. First, a similar argument shows that
the fully Bayes posterior is not KL-convergent to the so-called
``oracle posterior'' $\pi(\bolds{\theta}\mid\mathbf{y}, \mu
_T)$---that is, the
conditional posterior distribution for $\bolds{\theta}$, given the
true value
of $\mu$. This is not a source of worry for Bayesians, but it makes
clear that the disagreement between EB and FB procedures cuts both
ways, and is not merely a ``failure'' of empirical-Bayes; if a
non-Bayesian's goal is to reconstruct
the oracle posterior, this could be achieved by empirical-Bayes
analysis but not by full Bayes.

Second, the situation described above has the sample size $n$ equal to
the number of unknown parameters $m$. If $n$ grows relative to $m$, the
full Bayes and empirical-Bayes/oracle posteriors can indeed be
KL-convergent. For instance, suppose there are $r$ independent
replicate observations for each $\mu_i$. Then a similar calculation
shows that $\operatorname{KL}(\pi_F \parallel\pi_E) = O ( 1/r )$ as $ r
\rightarrow\infty$, so that KL convergence between the two approaches
would obtain.

\section{Results for variable selection}
\label{KLforlinearmodels2}

For the variable-selection problem, explicit expressions for the KL
divergence between empirical-Bayes and fully Bayes procedures are not
available. It is also quite difficult to characterize the sampling
distribution of $\hat{p}$, the empirical-Bayes estimate for the prior
inclusion probability $p$. It is therefore not yet possible to give a
general characterization of whether, and when, the empirical-Bayes
variable-selection procedure is KL-convergent, in the sense defined
above, to a fully Bayesian procedure.

Three interesting sets of results are available, however. First and
most simply, we can characterize the KL divergence between the
prior probability distributions of the fully Bayesian and
empirical-Bayesian procedures. Second, we can characterize the limiting
\textit{expected} Kullback--Leibler divergence between EB and FB
posteriors, even if we cannot characterize the limiting KL divergence
itself. Third, we can compare the asymptotic behavior of the full Bayes
and empirical-Bayes prior model probabilities for models in a size
neighborhood of the true model.

We denote the empirical-Bayes prior distribution over model indicators
by $p_E(M_{\bolds{\gamma}})$ and the fully-Bayesian distribution
(with uniform prior on $p$) by $p_F(M_{\bolds{\gamma}})$.
Similarly, after observing data $D$, we write $p_E(M_{\bolds{\gamma
}} \mid\mathbf{Y})$ and $p_F(M_{\bolds{\gamma}} \mid\mathbf{Y})$
for the posterior distributions.

\subsection{Prior KL divergence}
The first two theorems prove the existence of lower bounds on how close
the EB and FB priors can be, and show that these lower bounds become
arbitrarily large as the number of tests $m$ goes to infinity. We refer
to these lower bounds as ``information gaps,'' and give them in both
Kullback--Leibler (Theorem \ref{lowerboundKL}) and Hellinger (Theorem
\ref{lowerboundHellinger}) flavors.
\begin{theorem} \label{lowerboundKL}
Let $\underline{\mathrm{G}}(m) = \min_{\hat{p}} \operatorname
{KL}(p_F(M_{\bolds
{\gamma}}) \parallel p_E(M_{\bolds{\gamma}}) )$. Then
$\underline{\mathrm{G}}(m) \to\infty$ as $m \to\infty$.
\end{theorem}
\begin{pf}
The KL divergence is
%
%
\begin{eqnarray}
\label{priorKL}
\operatorname{KL} &=& \sum_{k=0}^m \frac{1}{m+1} \biggl[ \log\biggl( \frac{1}{m+1}
\pmatrix{m \cr k}^{-1} \biggr)
- \log\bigl( \hat{p}^{k} \cdot(1 - \hat{p})^{m - k} \bigr) \biggr]
\nonumber\\
&=& - \log(m+1) \\
&&{}- \frac{1}{m+1} \sum_{k=0}^m \biggl[ \log\pmatrix{m \cr k}
+ k \log\hat{p} + (m-k) \log(1 - \hat{p}) \biggr] .\nonumber
\end{eqnarray}
This is minimized for $\hat{p} = 1/2$ regardless of $m$, meaning that
%
%
\begin{eqnarray} \label{minKLgap}
\underline{\mathrm{G}}(m) &=& - \log(m+1) - \frac{1}{m+1}
\sum_{k=0}^m \biggl[ \log\pmatrix{m \cr k} + m \log(1/2) \biggr]
\nonumber\\[-8pt]\\[-8pt]
&=& m \log2 - \log(m+1) - \frac{1}{m+1} \sum_{k=0}^m
\log\pmatrix{m \cr k} .\nonumber
\end{eqnarray}
The first (linear) term in (\ref{minKLgap}) dominates the second
(logarithmic) term, whereas results in \citet{gould1964} show the
third term to be asymptotically linear in $m$ with slope $1/2$. Hence,
$\underline{\mathrm{G}}(m)$ grows linearly with $m$, with asymptotic positive
slope of $\log2 - 1/2$.
\end{pf}
\begin{theorem} \label{lowerboundHellinger}
Let $\underline{\mathrm{H}}^2(m) = \min_{\hat{p}}
\mathrm{H}^2(p_F(M_{\bolds{\gamma}}) \parallel p_E(M_{\bolds{\gamma
}}) )$. Then $\underline{\mathrm{H}}^2(m) \to1$ as $m \to
\infty$.
\end{theorem}
\begin{pf}
%
%
\begin{equation}
\mathrm{H}^2(p_F(M_{\bolds{\gamma}}) \parallel p_E(M_{\bolds
{\gamma}})) = 1 - \frac{1}{\sqrt{m+1}} \sum_{k=0}^m \sqrt{ \pmatrix{m
\cr
k} \hat{p}^k (1- \hat{p})^{m-k} } .
\end{equation}
This distance is also minimized for $\hat{p} = 1/2$, meaning that
%
%
\begin{equation}
\underline{\mathrm{H}}^2(m)
= 1 - (m+1)^{-1/2} \cdot2^{-m/2} \cdot\sum_{k=0}^m \sqrt{ \pmatrix{m
\cr
k} } .
\end{equation}
A straightforward application of Stirling's approximation to the
factorial function shows that
%
%
\begin{equation}
\lim_{m \to\infty} \Biggl[ (m+1)^{-1/2} \cdot2^{-m/2} \cdot\sum
_{k=0}^m \sqrt{ \pmatrix{m \cr k}} \Biggr] = 0 ,
\end{equation}
from which the result follows immediately.
\end{pf}

In summary, the ex-post prior distribution associated with the
EB procedure is particularly troubling when the number of tests $m$
grows without bound. On the one hand, when the true value of $k$
remains fixed or grows at a rate slower than $m$---that is, when
concerns over false positives become the most trenchant, and the case
for a Bayesian procedure exhibiting strong multiplicity control becomes
the most convincing---then $\hat{p} \to0$ and the EB prior
$p_E(M_{\bolds{\gamma}})$ becomes arbitrarily bad as an
approximation to $p_F(M_{\bolds{\gamma}})$. (Here, the correction
under the empirical-Bayes approach will be more aggressive compared
with the Bayesian approach, and some may consider this additional
aggressiveness to be a source of strength.) On the other hand, if the
true $k$ is growing at the same rate as $m$, then the best one can hope
for is that $\hat{p} = 1/2$. And even then, the information gap between
$p_F(M_{\bolds{\gamma}})$ and $p_E(M_{\bolds{\gamma}})$ grows
linearly without bound (for KL divergence), or converges to 1 (for
Hellinger distance).

\subsection{Posterior KL divergence}

We now prove a theorem showing that, under very mild conditions, the
expected KL divergence between FB and EB posteriors for the
variable-selection problem is infinite. This version assumes that the
error precision $\phi$ is fixed, but the generalization to an unknown
$\phi$ is straightforward.
\begin{theorem} \label{unboundedposteriorKL}
In the variable-selection problem, let $m$, $n > m$, and $\phi> 0$ be
fixed. Suppose $\mathbf{X}_{\bolds{\gamma}}$ is of full rank for
all models and that the family of priors for model-specific parameters,
$ \{ \pi(\bolds{\beta}_{\bolds{\gamma}}) \}$, is such that
$p(\bolds{\beta}_{\bolds{\gamma}} = \mathbf{0}) < 1$ for all
$M_{\bolds{\gamma}}$. Then, for any true model $M_{\bolds
{\gamma}}^T$, the expected posterior KL divergence $\mathrm{E}[
\operatorname{KL}\{p_F(M_{\bolds{\gamma}} \mid\mathbf{Y})
\parallel
p_E(M_{\bolds{\gamma}} \mid\mathbf{Y}) \}]$ under this true model
is infinite.
\end{theorem}
\begin{pf}
The posterior KL divergence is
%
%
\begin{equation}\qquad
\operatorname{KL}(p_F(M_{\bolds{\gamma}} \mid\mathbf{Y}) \parallel
p_E(M_{\bolds{\gamma}} \mid\mathbf{Y})) = \sum_{\Gamma}
p_F(M_{\bolds{\gamma}} \mid\mathbf{Y}) \cdot\log
\biggl( \frac{p_F(M_{\bolds{\gamma}} \mid\mathbf{Y})}
{p_E(M_{\bolds{\gamma}} \mid\mathbf{Y})} \biggr) .
\end{equation}
This is clearly infinite if there exists a model $M_{\bolds{\gamma
}}$ for which $p_E(M_{\bolds{\gamma}} \mid\mathbf{Y}) = 0$ but
$p_F(M_{\bolds{\gamma}} \mid\mathbf{Y}) > 0$. Since the fully
Bayesian posterior assigns nonzero probability to all models, this
condition is met whenever the empirical-Bayesian solution is $\hat{p} =
0$ or $\hat{p} = 1$.
Thus, it suffices to show that $\hat{p}$ will be $0$ with positive
probability under any true model.

Assume without loss of generality that $\phi= 1$. Recall that we are
also assuming that $\pi(\alpha) = 1$ for all models, and that the
intercept is orthogonal to all other covariates. Letting $\bolds
{\beta}_{\bolds{\gamma^*}} = (\alpha, \bolds{\beta}_{\bolds
{\gamma}})^t$ for model
$M_{\bolds{\gamma}}$, and letting $L(\cdot)$ stand for the
likelihood, the marginal likelihood for any model can then be written
%
%
\begin{equation}
\label{marglikegivenphi}
f(\mathbf{Y} \mid M_{\bolds{\gamma}}) =
L(\hat{\bolds{\beta}}_{\gamma}^*)
\cdot
\sqrt{2\pi/n}
\int_{\mathbb{R}^{k_{\gamma}}}
g(\bolds{\beta}_{\bolds{\gamma}})\pi(\bolds{\beta}_{\bolds
{\gamma}}) \,\dd
\bolds{\beta}_{\bolds{\gamma}} ,
\end{equation}
where
\[
g(\bolds{\beta}_{\bolds{\gamma}})
= \exp\bigl\{
-\tfrac{1}{2}
(\bolds{\beta}_{\bolds{\gamma}} - \hat{\bolds{\beta}}_{\gamma} )^t
\mathbf{X}_{\bolds{\gamma}}^t \mathbf{X}_{\bolds{\gamma}}
(\bolds{\beta}_{\bolds{\gamma}} - \hat{\bolds{\beta}}_{\gamma} )\bigr\}.
\]


The Bayes factor for comparing the null model to any model is
\[
B_{\gamma}(\mathbf{Y}) =
\frac{f(\mathbf{Y} \mid M_0)}
{f(\mathbf{Y} \mid M_{\bolds{\gamma}}) } ,
\]
which from (\ref{marglikegivenphi}) is clearly continuous as a function
of $\mathbf{Y}$ for every $\bolds{\gamma}$. Evaluated at $\mathbf
{Y} = \mathbf{0}$ (so that $\hat{\bolds{\beta}}_{\gamma}$ then
equals $0$), this Bayes factor satisfies
%
%
\begin{equation}
B_{\gamma}( \mathbf{0})
=
\biggl(
\int_{\mathbb{R}^{k_{\gamma}}}
\exp\biggl\{
-\frac{1}{2}
\bolds{\beta}_{\bolds{\gamma}}^t
\mathbf{X}_{\bolds{\gamma}}^t \mathbf{X}_{\bolds{\gamma}}
\bolds{\beta}_{\bolds{\gamma}}
\biggr\}
\pi(\bolds{\beta}_{\bolds{\gamma}}) \,\dd\bolds{\beta}_{\bolds
{\gamma}}
\biggr)^{-1} > 1
\end{equation}
for each $M_{\bolds{\gamma}}$ under the assumptions of the theorem.

By continuity, for every model $M_{\bolds{\gamma}}$ there exists
an $\varepsilon_{\gamma}$ such that $B_{\gamma}(\mathbf{Y}) > 1$
for any
$|\mathbf{Y}| < \varepsilon_{\gamma}$. Let $\varepsilon^* = \min
_{\bolds
{\gamma}} \varepsilon_{\gamma}$. Then for $\mathbf{Y}$ satisyfing
$|\mathbf
{Y}| < \varepsilon^*$, $B_{\gamma}(\mathbf{Y}) > 1$ for all nonnull
models, meaning that $M_0$ will have the largest marginal likelihood.
By Lemma \ref{EBextremelemma}, $\hat{p} = 0$ when such a $\mathbf{Y}$
is observed.

But under any model, there is positive probability of observing
$|\mathbf{Y}| < \varepsilon^*$ for any positive $\varepsilon^*$,
since this
set has positive Lebesgue measure. Hence, regardless of the true model,
there is positive probability that the KL divergence
$\operatorname{KL}(p_F(M_{\bolds{\gamma}} \mid\mathbf{Y})
\parallel
p_E(M_{\bolds{\gamma}} \mid\mathbf{Y}) )$ is infinite under the
sampling distribution $p(\mathbf{Y} \mid M_{\bolds{\gamma}})$, and
so its expectation is clearly infinite.
\end{pf}

Since the expected KL divergence is infinite for any number $m$ of
variables being tested, and for any true model, it is clear that
$\mathrm{E}(\operatorname{KL})$ does not converge to $0$ as $m \to\infty$. This, of
course, is a weaker conclusion than would be a lack of KL convergence
in probability.

In Theorem \ref{unboundedposteriorKL} the expectation is with respect
to the sampling distribution under a specific model $M_{\bolds
{\gamma}}$, with $\bolds{\beta}_{\bolds{\gamma}}$ either
fixed or marginalized away with respect to a prior distribution. But
this result implies an infinite expectation with respect to other
reasonable choices of the expectation distribution---for example, under
the Bernoulli sampling model for $\bolds{\gamma}$ in (\ref
{conditionalmodelprior}) with fixed prior inclusion probability $p$.

\subsection{Asymptotic behavior of prior model probabilities under EB
and FB procedures}

While interesting, the results in the previous two sections do not
consider the usual type of asymptotic comparison, namely, how do the
full Bayes and empirical Bayes posterior distributions converge as $m
\rightarrow\infty$? It is not clear that such asymptotic comparisons
are possible in general, although very interesting results can be
obtained in particular contexts [cf. \citet{bogdanghosh2008b},
\citet{bogdanetal2008}].

A rather general insight related to such comparison can be obtained,
however, by focusing on the prior probabilities
of ``high posterior'' models, as $m \rightarrow\infty$. To do so, we
first need an approximation to the full Bayes prior probability of
$M_{\bolds{\gamma}}$, given in the following lemma. The proof is
straightforward Laplace approximation, and is omitted.
\begin{lemma}
\label{lemma.laplace-approx}
As $m \rightarrow\infty$, consider models of size $k_{\bolds
{\gamma}}$ such that $k_{\bolds{\gamma}}/m$ is bounded away from 0
and 1. Then the Bayesian prior probability of $M_{\bolds{\gamma}}$
with prior $\pi(p)$ is
\begin{eqnarray*}
p_F(M_{\bolds{\gamma}}) &=& \int_0^1 p(M_{\bolds{\gamma}}
\mid p) \pi(p) \,\dd p \\
&=&
\biggl(\frac{k_{\bolds{\gamma}}}{m} \biggr)^{k_{\bolds{\gamma}}}
\biggl(1- \frac{k_{\bolds{\gamma}}}{m} \biggr)^{m-k_{\bolds
{\gamma}}}
\biggl[\frac{(2\pi)({k_{\bolds{\gamma}}}/{m})
(1-{k_{\bolds{\gamma}}}/{m})\pi({k_{\bolds{\gamma
}}}/{m})}{m} \biggr]^{1/2}\\
&&{}\times
\{1+o(1)\} ,
\end{eqnarray*}
providing $\pi(\cdot)$ is continuous and nonzero.
\end{lemma}

Now suppose $p_T$ is the true prior variable inclusion probability and
consider the most favorable situation for empirical Bayes analysis, in
which the empirical Bayes estimate for $p_T$ satisfies
%
%
\begin{equation}
\label{eq.EB-optimal}
\hat{p} = p_T (1+\varepsilon_{E})\qquad \mbox{where }\varepsilon_{E}\mbox{ is
} O \biggl(\frac{1}{\sqrt{m}} \biggr)\mbox{ as } m \rightarrow\infty.
\end{equation}
It is not known in general when this holds, but it does hold in
exchangeable contexts where each variable is in or out of the model
with unknown probability $p$, since such problems are equivalent to
mixture model problems.

For models far from the true model, the prior model probabilities given
by the Bayesian and empirical Bayesian approaches can be extremely
different. Hence, it is most interesting to focus on models that are
close to the true model for the comparison. In particular, we restrict
attention to models whose size differs from the true model by $O(\sqrt{m})$.
\begin{theorem}
Suppose the true model size $k_T$ satisfies $k_T/m = p_T +O(1/\sqrt
{m})$ as $m \rightarrow
\infty$, where $0< p_T <1$. Consider all models $M_{\bolds{\gamma
}}$ such that $k_T-k_{\bolds{\gamma}} =O(\sqrt{m})$, and consider
the optimal situation for EB in which (\ref{eq.EB-optimal}) holds. Then
the ratio of the prior probabilities assigned to such models by the
Bayes approach and the empirical Bayes approach satisfies
\begin{eqnarray*}
\frac{p_F(M_{\bolds{\gamma}})}{p_E(M_{\bolds{\gamma}})}
&=& \biggl(\frac{k_{\bolds{\gamma}}}{m} \biggr)^{k_{\bolds
{\gamma}}}
\biggl(1- \frac{k_{\bolds{\gamma}}}{m} \biggr)^{m-k_{\bolds
{\gamma}}}
\biggl[(2\pi)\biggl(\frac{k_{\bolds{\gamma}}}{m}\biggr)
\biggl(1-\frac{k_{\bolds{\gamma}}}{m}\biggr)\pi\biggl(\frac{k_{\bolds{\gamma
}}}{m}\biggr) \biggr]^{1/2}\\
&&{}\times m^{-1/2} \{1+o(1)\} \\
&&{}\times\bigl(
{ (\hat{p})^{k_{\bolds{\gamma}}}(1- \hat{p})^{m-k_{\bolds
{\gamma}}}}\bigr)^{-1} \\
&=& O \biggl(\frac{1}{\sqrt{m}} \biggr) ,
\end{eqnarray*}
providing $\pi(\cdot)$ is continuous and nonzero.
\end{theorem}
\begin{pf}
Note that
%
%
\begin{equation}\qquad
p_{E}(M_{\bolds{\gamma}})= (\hat{p})^{k_{\bolds{\gamma}}}(1-
\hat{p})^{m-k_{\bolds{\gamma}}}
= \{p_T (1+\varepsilon_{E})\}^{k_{\bolds{\gamma}}} \{1- p_T
(1+\varepsilon_{E})\}^{m-k_{\bolds{\gamma}}}.
\end{equation}
Taking the $\log$ and performing a Taylor expansion yields
\begin{eqnarray*}
\log{p_{E}(M_{\bolds{\gamma}})} &=&
\log\{ p_T^{k_{\bolds{\gamma}}}(1- p_T)^{m-k_{\bolds
{\gamma}}} \} + {k_{\bolds{\gamma}}} \log{(1+\varepsilon_{E})}
\\
&&{}+ {(m-k_{\bolds{\gamma}}}) \log\biggl\{1- \frac{p_T}{(1-p_T)}
\varepsilon_{E} \biggr\} \nonumber\\
&=& \log\{ p_T^{k_{\bolds{\gamma}}}(1-
p_T)^{m-k_{\bolds{\gamma}}} \} + {k_{\bolds{\gamma}}}
\{ \varepsilon_{E} + O(\varepsilon_{E}^2)\} \\
&&{}+(m-{k_{\bolds{\gamma}}}) \biggl\{ - \frac{p_T}{(1-p_T)} \varepsilon_{E}
+ O(\varepsilon_{E}^2) \biggr\} \nonumber\\
&=& \log\{ p_T^{k_{\bolds{\gamma}}}(1-
p_T)^{m-k_{\bolds{\gamma}}} \} + \bigl\{ k_T + O(\sqrt{m}) \bigr\} \{
\varepsilon_{E} + O(\varepsilon_{E}^2) \}
\nonumber\\
&&{}+ \bigl\{m - k_T - O(\sqrt{m})\bigr\} \biggl\{ - \frac{p_T}{(1-p_T)} \varepsilon_{E}
+ O(\varepsilon_{E}^2) \biggr\} \nonumber\\
&=& \log\{ p_T^{k_{\bolds{\gamma}}}(1-
p_T)^{m-k_{\bolds{\gamma}}} \} + O\bigl(\sqrt{m} \varepsilon_{E}\bigr) +
O(m\varepsilon_{E}^2) \nonumber\\
&=& \log\{ p_T^{k_{\bolds{\gamma}}}(1-
p_T)^{m-k_{\bolds{\gamma}}} \} + O(1) .
\end{eqnarray*}

A nearly identical argument using Lemma \ref{lemma.laplace-approx}
shows that the $\log$ Bayesian prior probability for these models is
%
%
\begin{equation}
\log\{ p_F(M_{\bolds{\gamma}}) \} = \log\{ p_T^{k_{\bolds
{\gamma}}}(1- p_T)^{m-k_{\bolds{\gamma}}} \} - \log{\sqrt{m}} +
O(1) ,
\end{equation}
from which the result is immediate.
\end{pf}

So we see that, even under the most favorable situation for the
empirical-Bayes analysis, and even when only considering models that
are close to the true model in terms of model size, the prior
probabilities assigned by the Bayes approach are smaller by a factor of
order $1/\sqrt{m}$ than those assigned by the empirical-Bayes approach.
The effect of this very significant difference in prior probabilities
will be context dependent,
but the result does provide a clear warning that the full Bayes and
empirical-Bayes answers can differ---even when $m \rightarrow\infty$
and even when there is sufficient information in the data to guarantee
the existence of a consistent estimator for $p_T$.

The theorem also shows that the empirical-Bayes procedure provides a
better asymptotic approximation to the ``oracle'' prior probabilities,
which may be argued by some to be the main goal of empirical-Bayes
analysis. At least for this ideal scenario, the EB approach assigns
larger prior probabilities to models which are closer to the true
model. Of course, this fact is not especially relevant from the fully
Bayesian perspective, and does not necessarily counterbalance the
problems associated with ignoring uncertainty in the estimator for $p$.

Finally, this difference in prior probabilities will not always have a
large effect. For instance, if $n \rightarrow\infty$ at a fast enough
rate compared with $m$, then the Bayes and empirical-Bayes approach
will typically agree simply because all of the posterior mass will
concentrate on a single model [i.e., one of the marginal likelihoods
$f(\mathbf{Y} \mid M_{\bolds{\gamma}})$ will become dominant], and
so the assigned prior probabilities will be irrelevant.

\section{Numerical investigation of empirical-Bayes variable selection}
\label{numerical-studies}

This section presents numerical results that demonstrate practical,
finite-sample significance of some of the qualitative differences
mentioned above. As in the previous section, most of the investigation
is phrased as a comparison of empirical-Bayes and full Bayes, taken
from the fully Bayesian perspective.

Note that $m$ here is taken to be moderate (14 for the simulation study
and 22 for the real data set); the intent is to focus on the magnitude
of the difference that one can expect in variable selection problems of
such typical magnitude. Of course, such $m$ are not large enough that
one would automatically expect the empirical-Bayes approach to provide
an accurate estimate of $p$, and so differences are to be expected, but
it is still useful to see the magnitude of the differences. For a
larger $m$ situation, see Table \ref{postincfakedata}; for the largest $m$ in that table,
the full Bayes and empirical-Bayes answers are much closer. The
rationale for taking these values of $m$ is that they allow the model
space to be enumerated, avoiding potential confounding effects due to
computational difficulties.

\subsection{Results under properly specified priors}
\label{ebgoeswrong}


The following simulation was performed 75,000 times for each of four
different sample sizes:

\begin{enumerate}
\item Draw a random $m \times n$ design matrix $\mathbf{X}$ of
independent $\mathrm{N}(0,1)$ covariates.
\item Draw a random $p \sim\mathrm{U}(0,1)$, and draw a sequence of $m$
independent Bernoulli trials with success probability $p$ to yield a
binary vector $\bolds{\gamma}$ encoding the true set of regressors.
\item Draw $\bolds{\beta}_{\bolds{\gamma}}$, the vector of
regression coefficients corresponding to the nonzero elements of
$\bolds{\gamma}$, from a Zellner--Siow prior. Set the other
coefficients $\bolds{\beta}_{-\bolds{\gamma}}$ to $0$.
\item Draw a random vector of responses $\mathbf{Y} \sim\mathrm
{N}(\mathbf{X}
\bolds{\beta}, \mathbf{I})$.
\item Using only $\mathbf{X}$ and $\mathbf{Y}$, compute marginal
likelihoods (assuming Zellner--Siow priors) for all $2^m$ possible
models; use these quantities to compute $\hat{p}$ along with the EB and
FB posterior distributions across model space.
\end{enumerate}

In all cases $m$ was fixed at $14$, yielding a model space of size
$16\mbox{,}384$---large enough to be interesting, yet small enough to be
enumerated 75,000 times in a row. We repeated the experiment for four
different sample sizes ($n = 16$, $n=30$, $n=60$ and $n=120$) to
simulate a variety of different $m/n$ ratios.


Two broad patterns emerged from these experiments.

\begin{figure}

\includegraphics{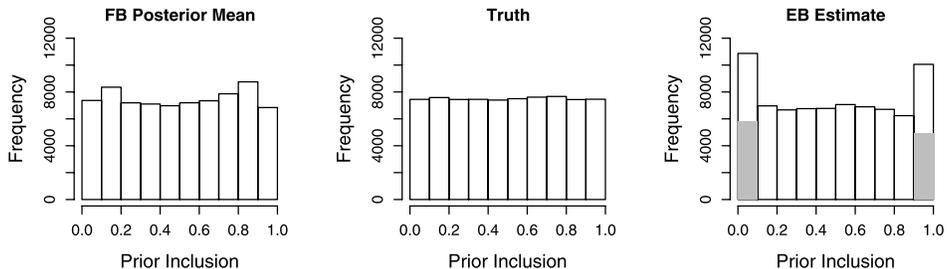}

\caption{Distribution of $\hat{p}$ in
the simulation study ($n=60$) with a correctly specified (uniform)
prior for $p$. The gray bars indicated the number of times, among
values of $\hat{p}$ in the extremal bins, that the empirical-Bayes
solution collapsed to the degenerate
$\hat{p} = 0$ or $\hat{p} = 1$.}\label{rhatp14}
\end{figure}

First, as Figure \ref{rhatp14} shows, the EB procedure gives the
degenerate $\hat{p} = 0$ or $\hat{p} = 1$ solution much too often. When
$n=60$, for example, almost $15 \%$ of cases collapsed to $\hat{p} = 0$
or $\hat{p} = 1$. This is essentially the same fraction of degenerate
cases as when $n=16$, which was $16\%$. This suggests that the issues
raised by Theorem~\ref{unboundedposteriorKL} can be quite serious in
practice, even when $n$ is large compared to $m$.

%
%

%
%
\begin{figure}

\includegraphics{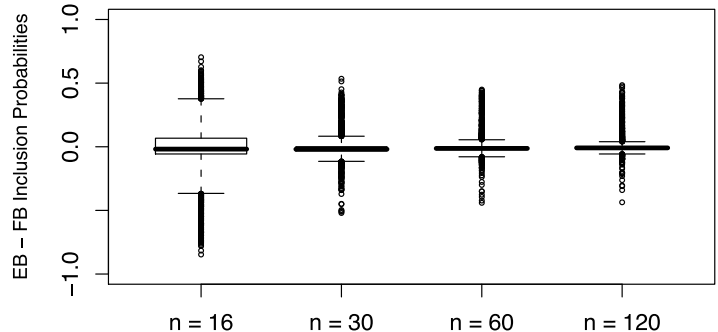}

\caption{Differences in all $m$
inclusion probabilities between EB and FB analyses across all
nondegenerate cases (i.e., where the EB solution does not collapse
to the boundary). The percentage of points lying outside the boxplot
whiskers ($1.5$ times the inter-quartile range) are as follows: $14\%$
for $n=16$, $12\%$ for $n=30$, $8\%$ for $n=60$ and $7\%$ for
$n=120$.}\label{incdiffp14}
\end{figure}

Second, even in nondegenerate situations, the two procedures often
reached very different conclusions about which covariates were
important. Figure \ref{incdiffp14} shows frequent large discrepancies
between the posterior inclusion probabilities given by the EB and FB
procedures. This happened even when $n$ was relatively large compared
to the number of parameters being tested, suggesting that even large
sample sizes do not render a data set immune to this difference. (Note
that Figure \ref{incdiffp14} only depicts the differences that arise
when the empirical-Bayes solution does not collapse to either $0$ or $1$.)




\subsection{Results under improperly specified priors}
\label{wrong-priors}

The previous section\break demonstrated that significant differences can
exist between fully Bayesian and empirical-Bayes variable selection in
finite-sample settings. There was an obvious bias, however, in that the
fully Bayesian procedure was being evaluated under its true prior
distribution, with respect to which it is necessarily optimal.

%
%
\begin{figure}

\includegraphics{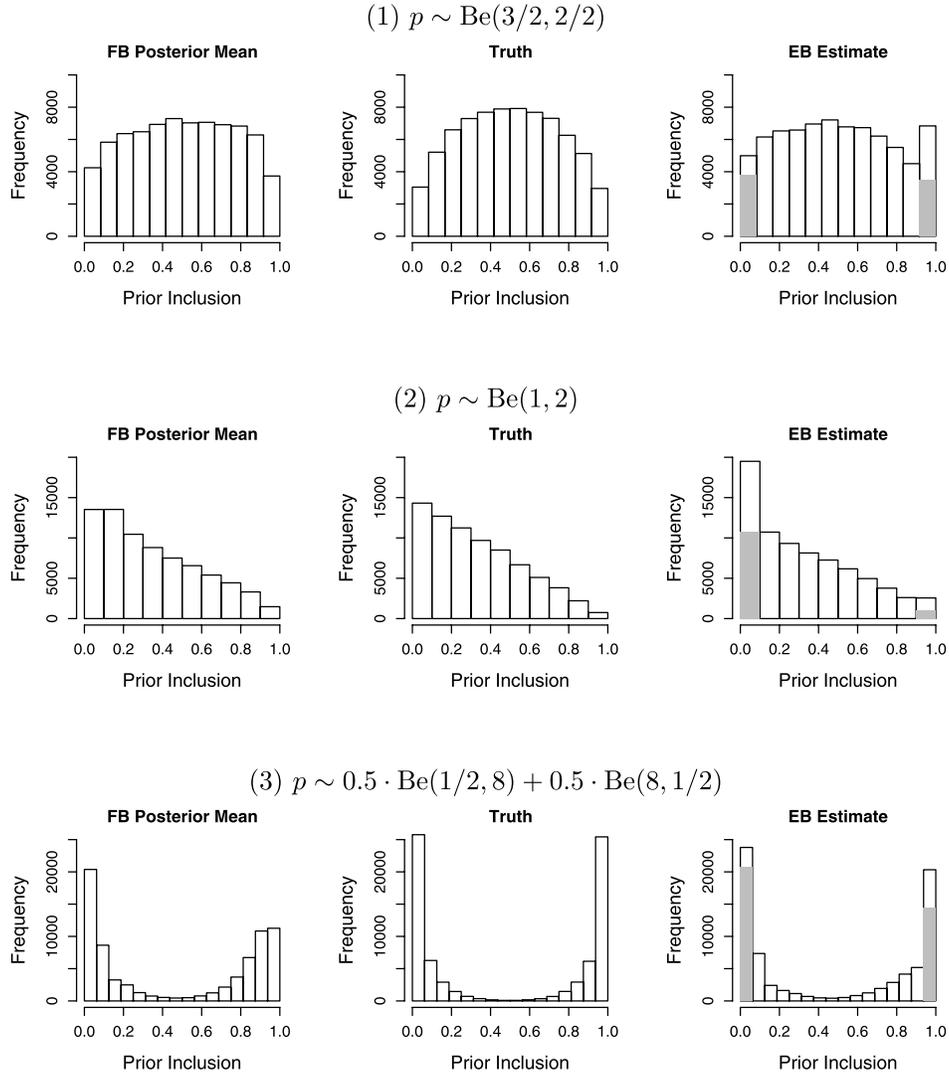}

\caption{Distribution of
$\hat{p}$ ($n=60$) in different versions of the simulation study, where
the fully Bayesian model had a misspecified (uniform) prior on $p$.
The gray bars indicated the number of times, among values of $\hat{p}$
in the extremal bins, that the empirical-Bayes solution collapsed to
the degenerate $\hat{p} = 0$ or $\hat{p} = 1$.}\vspace*{-4pt}\label{rhatp14wrongmodel}
\end{figure}

It is thus of interest to do a similar comparison for situations in
which the prior distribution is specified incorrectly: the fully
Bayesian answers will assume a uniform prior $p$, but $p$ will actually
be drawn from a nonuniform distribution. We limit ourselves to
discussion of the analogue of Figure \ref{rhatp14} for various
situations, all with $m=14$ and $n=60$. Three different choices of the
true distribution for $p$ were investigated, again with 75,000
simulated data sets each:
\begin{enumerate}
\item$p \sim\operatorname{Be}(3/2,3/2)$, yielding mainly moderate
(but not uniform)
values of $p$.
\item$p \sim\operatorname{Be}(1,2)$, yielding mainly smaller values
of $p$.
\item$p \sim0.5 \cdot\operatorname{Be}(1/2,8) + 0.5 \cdot
\operatorname{Be}(8,1/2)$, yielding
primarily values of $p$ close to 0 or 1.
\end{enumerate}

The results are summarized in Figure \ref{rhatp14wrongmodel}. In each
case the central pane shows the true distribution of $p$, with the left
pane showing the Bayesian posterior means under the uniform prior and
the right pane showing the empirical-Bayes estimates~$\hat{p}$.

As expected, the incorrectly specified Bayesian model tends to shrink
the estimated values of $p$ back to the prior mean of $0.5$. This
tendency is especially noticeable in Case 3, where the true
distribution contains many extreme values of~$p$. This gives the
illusion that empirical-Bayes tends to do better here.

Notice, however, the gray bars in the right-most panes. These bars
indicate the percentage of time, among values of $\hat{p}$ that fall in
the left- or right-most bins of the histogram, that the empirical-Bayes
solution is exactly $0$ or $1$, respectively. For example, of the
roughly 20,000 times that $\hat{p} \in[0, 0.1)$ in Case 2, it was
identically $0$ more than 10,000 of those times. (The fully Bayesian
posterior mean, of course, is never exactly 0 or 1.)

The bottom panel of Figure \ref{rhatp14wrongmodel} shows that,
paradoxically, where the fully Bayesian model is most incorrect, its
advantages over the empirical-Bayes procedure are the strongest. In the
mixture model giving many values of $p$ very close to $0$ or $1$,
empirical Bayes collapses to a degenerate solution nearly half the
time. Even if the extremal model is true in most of these cases, recall
that the empirical-Bayes procedure would result in an inappropriate
statement of certainty in the model. Of course, this would presumably
be noticed and some correction would be entertained, but the frequency
of having to make the correction is itself worrisome.

In these cases, while the fully Bayesian posterior mean is necessarily
shrunk back to the prior mean, this shrinkage is not very severe, and
the uniform prior giving rise to such shrinkage can easily be modified
if it is believed to be wrong. And in cases where the uniform prior is
used incorrectly, a slight amount of unwanted shrinkage seems a small
price to pay for the preservation of real prior uncertainty.

\subsection{Results when $p$ is fixed}
\label{subsection.pfixed}

\begin{figure}

\includegraphics{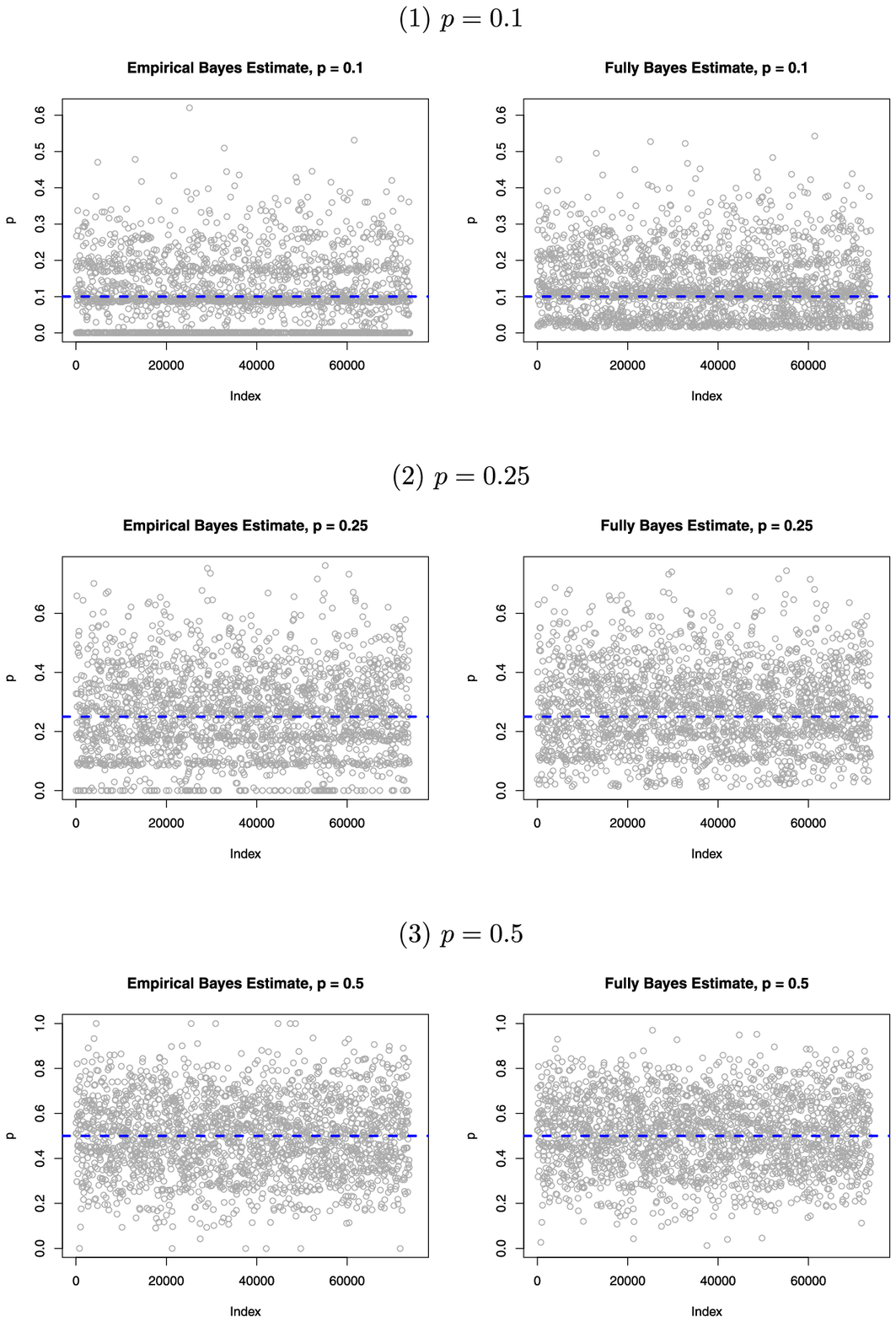}

\caption{Distribution of
$\hat{p}$ ($n=60$) in the fixed-$p$ versions of the simulation study
(2000 subsamples of the fake data sets). The dashed line indicates the
true value of $p$.}\label{rhatp14fixedp}
\end{figure}

We conducted a final version of the simulation with $p$ fixed at 3
different values: $p=0.10$, $p = 0.25$, and $p=0.5$. Figure \ref
{rhatp14fixedp} plots the estimated values of $p$ under the fully Bayes
and empirical-Bayes procedures. (For the sake of visual clarity only
the results from 2000 data sets are shown.)

%
%

It is clear that for the smallest value of $p=0.1$, the degenerate
solution $\hat{p} = 0$ occurs quite frequently. When $p$ is moderate
(as in the $0.25$ or $0.5$ cases), degeneracy occurs much less often.

It is also interesting to see the differences in how well the EB and FB
analysis approximate the ``oracle'' inclusion probabilities, which are
the posterior inclusion probabilities one would compute if one knew the
true Bernoulli probability $p$. This can be measured by looking at the
$\ell_1$ distance from the oracle estimate:
\[
\ell_1(\hat{\mathbf{p}}, \hat{\mathbf{p}}^{\mathrm{or}}) = \sum_{j=1}^{m}
| \hat
{p}_j - \hat{p}_j^{\mathrm{or}} | ,
\]
where $ \hat{p}_j^{\mathrm{or}} $ is the oracle posterior inclusion probability
for the $j$th variable.

The two procedures do quite similarly here, but with subtle
differences. For example, on the ``sparse'' ($p=0.1$) case, the mean
$\ell_1$ distance to the oracle answer across all Monte Carlo draws was
$0.36$ for the EB posterior, and $0.40$ for the FB posterior. Yet the
\textit{median} $\ell_1$ distance to the oracle answer was $0.27$ for the
FB posterior, and $0.30$ for the EB posterior.

These differences were largely consistent across other values of $p$.
This suggests that, while the FB procedure seems to reconstruct the
oracle posterior inclusion probabilities better for a larger number of
data sets (such as when the empirical-Bayes answer is degenerate), it
tends to miss by a larger amount than the EB procedure does. This
results in a worse level of average performance for the FB procedure in
reconstructing the oracle posterior inclusion probabilities.

\subsection{Example: Determinants of economic growth}
\label{econgrowthexample}

The following data set serves to illustrate the differences between EB
and FB answers in a scenario of typical size, complexity and $m/n$ ratio.

Many econometricians have applied Bayesian methods to the problem of
GDP-growth regressions, where long-term economic growth is explained in
terms of various political, social and geographical predictors.
\citet{FLS2001a} popularized the use of Bayesian model averaging
in the field; \citet{salaimartin2004} used a Bayes-like procedure
called BACE, similar to BIC-weighted OLS estimates, for selecting a
model; and \citet{leysteel2007} considered the effect of prior
assumptions (particularly the pseudo-objective $p = 1/2$ prior) on
these regressions.

We study a subset of the data from \citet{salaimartin2004}
containing 22 covariates on 30 different countries. A data set of this
size allows the model space to be enumerated and the EB estimate
$\hat{p}$ to be calculated explicitly, which would be impossible on the
full data set. The 22 covariates correspond to the top 10 covariates
flagged in the BACE study, along with 12 others chosen uniformly at
random from the remaining candidates.

Summaries of exact EB and FB analyses (with Zellner--Siow priors) can
be found in Table \ref{growthecon}. Two results are worth noting.
First, the EB inclusion probabilities are nontrivially different from
their FB counterparts, often disagreeing by $10\%$ or more.

%
%
\begin{table}[b]
\caption{Exact inclusion probabilities for 22 variables in a linear
model for GDP growth among\break a group of 30 countries}\label{growthecon}
\begin{tabular*}{\tablewidth}{@{\extracolsep{\fill}}lcc@{}}
\hline
\textbf{Covariate} & \textbf{Fully Bayes} & \textbf{Emp. Bayes} \\
\hline
East Asian dummy & 0.983 & 0.983 \\
Fraction of tropical area & 0.727 & 0.653 \\
Life expectancy in 1960 & 0.624 & 0.499 \\
Population density coastal in 1960s & 0.518 & 0.379 \\
GDP in 1960 (log) & 0.497 & 0.313 \\
Outward orientation & 0.417 & 0.318 \\
Fraction GDP in mining & 0.389 & 0.235 \\
Land area & 0.317 & 0.121 \\
Higher education 1960 & 0.297 & 0.148 \\
Investment price & 0.226 & 0.130 \\
Fraction confucian & 0.216 & 0.145 \\
Latin American dummy & 0.189 & 0.108 \\
Ethnolinguistic fractionalization & 0.188 & 0.117 \\
Political rights & 0.188 & 0.081 \\
Primary schooling in 1960 & 0.167 & 0.093 \\
Hydrocarbon deposits in 1993 & 0.165 & 0.093 \\
Fraction spent in war 1960--1990 & 0.164 & 0.095 \\
Defense spending share & 0.156 & 0.085 \\
Civil liberties & 0.154 & 0.075 \\
Average inflation 1960--1990 & 0.150 & 0.064 \\
Real exchange rate distortions & 0.146 & 0.071 \\
Interior density & 0.139 & 0.067 \\
\hline
\end{tabular*}
\end{table}

Second, if these are used for model selection, quite different results
would emerge. For instance, if median-probability models were selected
(i.e., one includes only those variables with inclusion probability
greater than $1/2$), the FB analysis would include the first four
variables (and would almost choose the fifth variable), while the EB
analysis would select only the first two variables (and almost the
third). While we would not endorse simply choosing a model here, note
that doing so would result in fundamentally different economic pictures
for the FB and EB analysis.


\section{Summary}
\label{summary}

This paper started out as an attempt to more fully understand when, and
how, multiplicity correction automatically occurs in Bayesian analysis,
and to examine the importance of ensuring that such multiplicity
correction is included. That the correction can only happen through the
choice of appropriate prior probabilities of models seemed to conflict
with the intuition that multiplicity correction occurs through
data-based adaptation of the prior-inclusion probability~$p$.

The resolution to this conflict---that the multiplicity correction is
indeed pre-fixed in the prior probabilities, but the amount of
correction employed will depend on the data---led to another conflict:
how can the empirical-Bayes approach to variable selection be an
accurate approximation to the full Bayesian analysis? Indeed, we have
seen in the paper that empirical-Bayes variable selection can lead to
results quite different than those from the full Bayesian analysis.
This difference was evidenced through examples (both simple pedagogical
examples and a more realistic practical example), through simulation
studies, and through information-based theoretical results. These
studies, as well as the results about the tendency of empirical-Bayes
variable selection to choose extreme $\hat p$, all supported the
general conclusions about empirical-Bayes variable selection that were
mentioned in the \hyperref[intro]{Introduction}.

\begin{appendix}\label{parameterpriors}

\section*{Appendix: Variations on Zellner's $g$-prior}

Conventional variable-selection priors rely upon the conjugate
normal-gamma family of distributions, which yields closed-form
expression for the marginal likelihoods. To give an appropriate scale
for the normal prior describing the regression coefficients,
\citet{zellner86} suggested a particular form of this family:
\begin{eqnarray*}
(\bolds{\beta} \mid\phi) &\sim& \mathrm{N}\biggl(\bolds{\beta}_0,
\frac{g}{\phi} (\mathbf{X}'\mathbf{X})^{-1} \biggr)  ,\\
\phi&\sim& \operatorname{Ga}\biggl(\frac{\nu}{2}, \frac{\nu s}{2} \biggr)
\end{eqnarray*}
with prior mean $\bolds{\beta}_0$, often chosen to be $0$. The
conventional choice $g = n$ gives a prior covariance matrix for the
regression parameters equal to the unit Fisher information matrix for
the observed data $\mathbf{X}$. This prior can be interpreted as
encapsulating the information arising from a single observation under a
hypothetical experiment with the same design as the one to be analyzed.

Zellner's $g$-prior was originally formulated for testing a precise
null hypothesis, $H_0\dvtx\bolds{\beta} = \bolds{\beta}_0$,
versus the alternative, $H_A\dvtx\bolds{\beta} \in\mathbb{R}^p$.
But others have adapted Zellner's methodology to the more general
problem of testing nested regression models by placing a flat prior on
the parameters shared by the two models and using a $g$-prior only on
the parameters not shared by the smaller model. This seems to run
afoul of the general injunction against improper priors in model
selection problems, but can nonetheless be formally justified by
arguments appealing to othogonality and group invariance; see, for
example, \citet{bergervarsh1998} and \citet{eaton1989}.
These arguments apply to cases where all covariates have been centered
to have a mean of zero, which is assumed without loss of generality to
be true.

A full variable-selection problem, of course, involves many nonnested
comparisons. Yet Bayes factors can still be formally defined using the
``encompassing model'' approach of \citet{zellnersiow1980}, who
operationally define all marginal likelihoods in terms of Bayes factors
with respect to a base model $M_B$:
%
%
\begin{equation}
\label{encompassingbayesfactor}
\operatorname{BF}(M_1 \dvtx M_2) = \frac{\operatorname{BF}(M_1 \dvtx
M_B)}{\operatorname{BF}(M_2 \dvtx
M_B)} .
\end{equation}

Since the set of common parameters which are to receive improper priors
depends upon the choice of base model, different choices yield a
different ensemble of Bayes factors and imply different ``operational''
marginal likelihoods. And while this choice of $M_B$ is free in
principle, there are only two such choices which yield a pair of nested
models in all comparisons: the null model and the full model.

In the null-based approach, each model is compared to the null model
consisting only of the intercept $\alpha$. This parameter, along with
the precision $\phi$, is common to all models, leading to a prior
specification that has become the most familiar version of Zellner's $g$-prior:
\begin{eqnarray*}
(\alpha, \phi\mid\bolds{\gamma}) &\propto& 1/\phi,\\
(\bolds{\beta}_{\bolds{\gamma}} \mid\phi, \bolds{\gamma
}) &\sim& \mathrm{N}\biggl(0, \frac{g}{\phi} (\mathbf{X}_{\bolds{\gamma
}}' \mathbf{X}_{\bolds{\gamma}})^{-1} \biggr) .
\end{eqnarray*}

This gives a simple expression for the Bayes factor for evaluating a
model $\bolds{\gamma}$ with $k$ regression parameters (excluding
the intercept):
%
%
\begin{equation}
\label{gnullbayesfactor}
\operatorname{BF}(M_{\bolds{\gamma}} \dvtx M_0) =
(1+g)^{(n-k_{\bolds
{\gamma}} - 1)/2} [1 + (1-R^2_{\bolds{\gamma}})g]^{-(n-1)/2} ,
\end{equation}
where $R^2_{\bolds{\gamma}} \in(0,1]$ is the usual coefficient of
determination for model $M_{\bolds{\gamma}}$.

Adherents of the full-based approach, on the other hand, compare all
models to the full model, on the grounds that the full model is usually
much more scientifically reasonable than the null model and provides a
more sensible yardstick [\citet{casellamoreno2002}]. This
comparison can
be done by writing the full model as
\[
M_F\dvtx\mathbf{Y} = \mathbf{X}^*_{\bolds{\gamma}} \theta_{\bolds
{\gamma}} + \mathbf{X}_{-\bolds{\gamma}} \bolds{\beta
}_{-\bolds{\gamma},}
\]
with the design matrix partitioned in the obvious way. Then a $g$-prior
is specified for the parameters in the full model not shared by the
smaller model, which again has $k$ regression parameters excluding the
intercept:
\begin{eqnarray*}
(\alpha, \bolds{\beta}_{\bolds{\gamma}}, \phi\mid\bolds
{\gamma}) &\propto& 1/\phi,\\
(\bolds{\beta}_{-\bolds{\gamma}} \mid\phi, \bolds{\gamma
}) &\sim& \mathrm{N}\biggl(0, \frac{g}{\phi} (\mathbf{X}_{-\bolds
{\gamma
}}' \mathbf{X}_{-\bolds{\gamma}})^{-1} \biggr).
\end{eqnarray*}

This does not lead to a coherent ``within-model'' prior specification
for the parameters of the full model, since their prior distribution
depends upon which submodel is considered. Nevertheless, marginal
likelihoods can still be consistently defined in the manner of
(\ref{encompassingbayesfactor}). Conditional upon $g$, this yields a
Bayes factor in favor of the full model of
%
%
\begin{equation}
\label{gfullbayesfactor}
\operatorname{BF}(M_F \dvtx M_{\bolds{\gamma}}) = (1 + g)^{(n - m -
1)/2} (1
+ gW)^{-(n - k - 1)/2},
\end{equation}
where $W = (1 - R_F^2)/(1 - R^2_{\bolds{\gamma}})$.

The existence of these simple expressions has made the use of
$g$-priors very popular. Yet $g$-priors yield display a disturbing type
of behavior often called the ``information paradox.'' This can be seen
in (\ref{gnullbayesfactor}): the Bayes factor in favor of
$M_{\bolds{\gamma}}$ goes to the finite constant $(1 + g)^{n - m -
1}$ as $R^2_{\bolds{\gamma}} \to1$ (which can only happen if
$M_{\bolds{\gamma}}$ is true and the residual variance goes to
$0$). For typical problems this will be an enormous number, but still
quite a bit smaller than infinity. Hence, the paradox: the Bayesian
procedure under a $g$-prior places an intrinsic limit upon the possible
degree of convincingness to be found in the data, a limit which is
confirmed neither by intuition nor by the behavior of the classical
test statistic.

\citet{liangpaulo07} detail several versions of
information-consistent $g$-like priors. One way is to estimate $g$ by
empirical-Bayes methods [\citet{georgefoster2000}]. A second, fully
Bayesian, approach involves placing a prior upon $g$ that satisfies the
condition $\int_0^{\infty} (1 + g)^{n - k_{\bolds{\gamma}} - 1}
\pi(g) \,\dd g = \infty$ for all $k_{\bolds{\gamma}} \leq p$,
which is a generalization of the condition given in
\citet{jeffreys1961} (see Chapter 5.2, equations 10 and 14).

This second approach generalizes the recommendations of
\citet{zellnersiow1980}, who compare models by placing a flat
prior upon common parameters and a $g$-like Cauchy prior on nonshared
parameters:
%
%
\begin{equation}
\label{zellnersiowpriors}
(\bolds{\beta}_{\bolds{\gamma}} \mid\phi) \sim C \biggl(0,
\frac{n}{\phi} (\mathbf{X}_{\bolds{\gamma}} ' \mathbf
{X}_{\bolds{\gamma}})^{-1} \biggr).
\end{equation}
These have come to be known as Zellner--Siow priors, and their use can
be shown to resolve the information paradox. Although they do not yield
closed-form expressions for marginal likelihoods, one can exploit the
scale-mixture-of-normals representation of the Cauchy distribution to
leave one-dimensional integrals over standard $g$-prior marginal
likelihoods with respect to an inverse-gamma prior, $g \sim
\operatorname{IG}(1/2,2/n)$. The Zellner--Siow null-based Bayes factor
under model
$M_{\bolds{\gamma}}$ then takes the form
%
%
\begin{eqnarray}
\label{ZSNmarglike}
\operatorname{BF}(M_{\bolds{\gamma}} \dvtx M_0) &=& \int_0^\infty
(1+g)^{(n-k_{\bolds{\gamma}} - 1)/2} [1 + (1-R^2_{\bolds
{\gamma}})g]^{-(n-1)/2}\nonumber\\[-8pt]\\[-8pt]
&&\hspace*{17.1pt}{}\times g^{-3/2} \exp\bigl(-n/(2g)\bigr) \,dg.\nonumber
\end{eqnarray}

A similar formula exists for the full-based version:
%
%
\begin{eqnarray}
\label{ZSFmarglike}
\operatorname{BF}(M_F \dvtx M_{\bolds{\gamma}}) &=& \int_0^\infty
(1+g)^{(n-m -
1)/2} [1 + Wg]^{-(n-k-1)/2}\nonumber\\[-8pt]\\[-8pt]
&&\hspace*{17.1pt}{}\times g^{-3/2} \exp\bigl(-n/(2g)\bigr) \,\dd g\nonumber
\end{eqnarray}
with $W$ given above.

These quantities can be computed by one-dimensional numerical
integration, but in high-dimensional model searches this will be a
bottleneck. Luckily there exists a closed-form approximation to these
integrals first noted in \citet{liangpaulo07}. It entails
computing the roots of a cubic equation, and extensive numerical
experiments show the approximation to be quite accurate. These Bayes
factors seem to offer an excellent compromise between good theoretical
behavior and computational tractability, thereby overcoming the single
biggest hurdle to the widespread practical use of Zellner--Siow priors.
\end{appendix}





%
\printaddresses

\mbox{}
\end{document}